\documentclass[12pt]{article}
\usepackage{amsfonts}

\usepackage{amsmath, amssymb}
\usepackage{amssymb}
\def\mineappendix{
        \setcounter{section}{1}
        \setcounter{subsection}{0}
        \def\thesection{\Alph{section}}
        \def\sectionap{\@startsection  {section}{1}{\z@}
                        {-3.5ex plus-1ex minus-.2ex} {0ex plus.2ex}
                        {\reset@font\Large\bf  Appendix:  \, }
                        }
        }
\makeatother
\def\Proclaim #1. #2\par{\bigbreak\noindent{\sc#1.\enspace}{\it#2}\par}

\newcommand{\gwi}[1]{\left< \, #1 \, \right>_{0}}
\newcommand{\gwii}[1]{\left< \hspace{-2pt} \left< \, #1 \,
        \right>  \hspace{-2pt} \right>_{0}}
\newcommand{\gwiis}[1]{\left< \hspace{-2pt} \left< \, #1 \,
        \right>  \hspace{-2pt} \right>_{0,\mathcal{H}}}

\newcommand{\gwione}[1]{\left< \, #1 \, \right>_{1}}
\newcommand{\gwiione}[1]{\left< \hspace{-2pt} \left< \, #1 \,
        \right> \hspace{-2pt} \right>_{1}}

\newcommand{\gwitwo}[1]{\left< \, #1 \, \right>_{2}}
\newcommand{\gwiitwo}[1]{\left< \hspace{-2pt} \left< \, #1 \,
        \right> \hspace{-2pt} \right>_{2}}

\newcommand{\gwig}[1]{\left< \, #1 \, \right>_{g}}
\newcommand{\gwih}[2]{ \left< \, #2 \, \right>_{#1}}
\newcommand{\gwiig}[1]{\left< \hspace{-2pt} \left< \, #1 \,
    \right> \hspace{-2pt} \right>_{g}}
\newcommand{\gwiih}[2]{\left< \hspace{-2pt} \left< \, #2 \,
    \right> \hspace{-2pt} \right>_{#1}}

\newcommand{\grav}[2]{\tau_{#1}(\gamma_{#2})}
\newcommand{\grava}[1]{\tau_{#1}(\gamma_{\alpha})}

\newcommand{\ga}{\gamma_{\alpha}}
\newcommand{\gua}{\gamma^{\alpha}}
\newcommand{\gb}{\gamma_{\beta}}
\newcommand{\gub}{\gamma^{\beta}}

\newcommand{\ve}{{\mathcal E}}
\newcommand{\vw}{{\cal W}}

\newtheorem{lem}{Lemma}[section]
\newtheorem{cor}[lem]{Corollary}
\newtheorem{thm}[lem]{Theorem}

\newtheorem{rem}[lem]{Remark}

\title{ Conditions for the vanishing of the genus-2 G-function  }

\author{Xiaobo Liu \thanks{Research of the first author was partially supported by NSF grant
DMS-0905227, NSFC Tianyuan special fund 11226027, SRFDP grant
20120001110051, and Peking University 985 fund.}, \quad  Xin Wang
\thanks{Research of the second author was partially supported by
SRFDP grant 20120001110051.}}

\date{}

\begin{document}
\maketitle
\begin{abstract}
In this paper we give some sufficient conditions for the vanishing of the genus-2 G-function, which was introduced by B. Dubrovin, S. Liu  and Y. Zhang in [DLZ]. As a corollary we prove their conjecture for the vanishing of the genus-2 G-function for ADE singularities.
\end{abstract}

In this paper, we study properties of the genus-2 potential function $F_{2}$ for semisimple Frobenius manifolds.
In \cite{DZ}, Dubrovin and Zhang proved that higher genus potential functions can be uniquely constructed from the Virasoro
constraints and the genus-0 potential function for semisimple Frobenius manifolds . They also gave an explicit formula for
$F_2$ for semisimple Frobenius manifolds under the assumption of the genus-2 Virasoro conjecture. In \cite{L07},
by studying genus-2 topological recursion relations for Gromov-Witten invariants, the first author of this paper
obtained a much simpler formula of $F_{2}$ in terms of genus-0 Gromov-Witten invariants for compact symplectic manifolds
with semisimple quantum cohomology, and proved the genus-2 Virasoro conjecture using this formula. Since the
proofs of \cite{L07} only used universal equations for Gromov-Witten invariants which was obtained from the
splitting principle and relations in the tautological rings of moduli spaces of stable curves, the formula
found in \cite{L07} holds for the genus-2 generating functions of descendant
invariants of any semisimple cohomological field theory. Both formulas in \cite{DZ} and \cite{L07} used canonical
coordinates for semisimple Frobenius manifolds (or a generalization of canonical coordinates for the big phase
space of the Frobenius manifolds). Since geometric invariants for cohomological field theories
usually appear as coefficients of the formal power series representations of potential functions in
terms of flat coordinates, it is desirable to have a formula for $F_2$ purely in terms of flat coordinates.

Recently in \cite{DLZ}, B. Dubrovin, S. Liu  and Y. Zhang decomposed $F_{2}$  as a summation of two functions:
 \begin{align}
 F_{2}=F'_{2}+G^{(2)}    \label{F2}
 \end{align}
where the first part $F'_{2}$ can be expressed purely in terms of flat coordinates, and the second part $G^{(2)}$, called the {\bf genus-$2$ G-function}, still depends on the canonical coordinates.
An explicit formula for the genus-2 G-function can be found in Appendix~\ref{sec:G2}. We will not define $F'_2$ since most part of this paper does not need this function. Interested readers are referred to \cite{DLZ} for the exact form of $F'_2$. It is an interesting problem to find conditions under which the genus-2 G-function vanishes, i.e.
$G^{(2)}=0$. In such cases, we will have a formula for $F_2$ completely expressed in terms of flat coordinates.

The main purpose of this paper is to give a set of geometric conditions which force the vanishing of the genus-2 G-function. Let $(\mathcal{H},\eta,\Lambda)$ be a cohomological field theory of rank $N$, where $\mathcal{H}$ is an N-dimensional vector space with a nondegenerate paring $\eta$ and $\Lambda:=\{\Lambda_{g,n}\}$ are n-linear $H^{\bullet}(\overline{\mathcal{M}}_{g,n})$-valued forms on
$\mathcal{H}$, i.e.
$
\Lambda_{g,n}\in Hom(\mathcal{H}^{\bigotimes n},H^{\bullet}(\overline{\mathcal{M}}_{g,n})),
$
 satisfying the axioms of cohomological field theory(cf: \cite{KM}).
Let $\{\gamma_{1},\gamma_{2},...,\gamma_{N}\}$ be a fixed basis of $\mathcal{H}$, and $\eta^{\alpha\beta}$ entries for the inverse of the matrix of $\eta$ in this basis. We will use $\eta^{\alpha \beta}$ to raise indices, for example $\gamma^{\alpha}=\sum_{\beta}\eta^{\alpha\beta}\gamma_{\beta}$ for any $\alpha$. Let
$\gwig{\grav{n_1}{\alpha_1} \cdots \grav{n_k}{\alpha_k}}$ be genus-$g$ descendant invariants of the cohomological field theory determined by $\Lambda$. We will identify $\grava{0}$ with $\gamma_{\alpha}$ and call
$\gwig{\gamma_{\alpha_1} \cdots \gamma_{\alpha_k}}$ primary invariants.
The genus-0 primary invariants of a cohomological field theory defines a Frobenius manifold structure on $\mathcal{H}$.
Consider the following three conditions:
\begin{align*}
&(C1) \quad\quad
\sum_{\alpha,\beta}\gwi{ {\gamma _\alpha }{\gamma ^\alpha }{\gamma _\beta }{\gamma ^\beta }\gamma_{\alpha_{1}}\cdot\cdot\cdot \gamma_{\alpha_{k}}} =0,
\\&(C2) \quad\quad  \gwione{ \gamma_{\alpha_{1}}\cdot\cdot\cdot \gamma_{\alpha_{k}}} =0,
\\&(C3) \quad\quad  \gwitwo{ \gamma_{\alpha_{1}}\cdot\cdot\cdot \gamma_{\alpha_{k}}} =0 \quad \mbox{and} \quad
\gwitwo{ \tau_{1}(\gamma_{\alpha_{1}})\gamma_{\alpha_{2}}\cdot\cdot\cdot \gamma_{\alpha_{k}}} =0,
\end{align*}
for all $\alpha_{1},\cdot\cdot\cdot,\alpha_{k}$, $k\geq 1$.
Condition (C1) is inspired by Lemma 2.5 in \cite{DLZ}. The relation between this condition and Lemma 2.5 in \cite{DLZ}
will be explained in Section~\ref{sec:pfcor}. The proof of Lemma 2.5 in \cite{DLZ} for $ADE$ singularities
was quite complicated and was done via case by case study (cf. Section 3.1 to Section 3.4 in \cite{DLZ}).
We will give a unified and much simpler proof of this lemma in Section~\ref{sec:pfcor}
(cf. Remark~\ref{rem:pfDLZ2.5}).
By the genus-1 constitutive relation of Dijkgraaf and Witten,
the genus-1 generating function $F_1$ can be written as a sum of a genus-0 function together with a function which still depends on genus-1 primary invariants (cf. \cite{DW}). The latter part of this decomposition was called the genus-1 G-function by Dubrovin and Zhang.
Condition (C2) is equivalent to that the genus-1 G-function is equal to a constant. The first part of the condition (C3) is a necessary condition for
the vanishing of genus-2 G-function since the restriction of $F'_2$ to the {\it small phase space} (i.e. $\mathcal{H}$) is zero.
Under the assumption of condition (C2), it can be shown that the second part of condition (C3)
is also a necessary condition for the vanishing of genus-2 G-function due to the string equation.
The main result of this paper is the following
\begin{thm} \label{thm:main}
For any semisimple cohomological field theory satisfying conditions (C1),(C2) and (C3),
 the genus-$2$ G-function vanishes.
\end{thm}

It can be shown that for the FJRW theory of $ADE$ singularities (cf. \cite{FJR}),
the three conditions (C1)-(C3) are all satisfied.
Therefore an immediate consequence of Theorem~\ref{thm:main} is the following
\begin{cor}  \label{cor:ADEG}
The genus-2 G-function vanishes for the cohomological field theories associated to $ADE$ singularities.
\end{cor}

It was conjectured in \cite{DLZ} that the genus-2 G-function vanishes for the cohomological field theories associated to
$ADE$ singularities and Gromov-Witten theory for $\mathbb{P}^{1}$-orbifolds of $ADE$ type. Corollary~\ref{cor:ADEG} gives
an affirmative answer to this
conjecture for $ADE$ singularities. Recently Y. Fu, S. Liu, Y. Zhang, and C. Zhou proved this conjecture for $A_n$
singularities using a completely different method (cf. \cite{FLZZ}). Their proof depends on specific structures of Frobenius
manifolds associated to $A_n$ singularities, which do not apply to singularities of $DE$ type.
In comparison, our method gives a more geometric explanation and applies to all $ADE$ singularities, and possibly also
applies to other cohomological field theories as well. We will discuss the case for $\mathbb{P}^{1}$-orbifolds of $ADE$
type in a forthcoming paper.

This paper is organized as follows.
In section~\ref{sec:prel}, we introduce notations and review basic theories needed in this paper.
We prove Theorem~\ref{thm:main} in section~\ref{sec:pfmain} and prove Corollary~\ref{cor:ADEG}
in section~\ref{sec:pfcor}.
In section~\ref{sec:rem}, we make some remarks on how to simplify $F_2$ if
the genus-2 G-function is zero and condition (C2) is satisfied.
In appendices, we give the precise definition for the genus-2 G-function,
the formula for genus-2 generating function as proved in \cite{L07},
precise formulas for some equations used in the proof of the main theorem, and two genus-1 3-point functions.

\section{Preliminaries}
\label{sec:prel}

Properties of the quantum product and idempotents on the big phase space for Gromov-Witten theory were studied in
\cite{L02}, \cite{L06} and \cite{L07}. These properties can be easily extended to any cohomological field theory.
In this section we will review concepts and properties needed for the proof of Theorem~\ref{thm:main}.

Let $(\mathcal{H},\eta,\Lambda)$ be a cohomological field theory of rank $N$. Fix a basis  $\{\ga \mid \alpha=1, \cdots N \}$ of $\mathcal{H}$ with $\gamma_{1}$ equal to the identity of the Frobenius manifold structure on $\mathcal{H}$.  As a convention, repeated Greek letter
indices are summed
up over their entire ranges from 1 to $N$.
The vector space $\mathcal{H}$ is also called {\it the small phase space}.
The {\it big phase space} is defined to be
$\mathcal{P}:=\prod_{n=0}^{\infty} \mathcal{H}$ and the corresponding
basis for the $n$-th copy of $\mathcal{H}$ in this product is denoted by
$\{ \grava{n} \mid \alpha=1, \ldots, N\}$ for $ n \geq 0$.
Let $t_{n}^{\alpha}$ be the coordinates on $\mathcal{P}$
with respect to the standard basis $\{ \grava{n} \mid \alpha=1, \ldots, N, \,\,\, n \geq 0\}$.
These coordinates are called the {\it flat coordinates} on
$\mathcal{P}$. We will also identify the small phase space $\mathcal{H}$ with
the subspace of $\mathcal{P}$ defined by $t_n^{\alpha}=0$ for $n>0$.

The genus-$g$ potential function $F_{g}$
 is
a formal power series of $t=(t_{n}^{\alpha})$ with coefficients given by
$\gwig{\grav{n_1}{\alpha_1} \cdots \grav{n_k}{\alpha_k}}$.
Derivatives of
$F_{g}$ with respect to $t_{n_{1}}^{\alpha_{1}}, \ldots, t_{n_{k}}^{\alpha_{k}}$ are denoted by
$\gwiig{\tau_{n_{1}}(\alpha_{1}) \cdots \tau_{n_{k}}(\alpha_{k})}$.
We will identify $\tau_{n}(\gamma_{\alpha})$ with $\frac{\partial}{\partial t_{n}^{\alpha}}$ as vector fields on the big phase space. We will also write $\tau_{0}(\gamma_{\alpha})$ simply as $\gamma_{\alpha}$.
Any vector field of the form $\sum_{\alpha}f_{\alpha}\gamma_{\alpha}$, where $f_{\alpha}$ are functions on the big phase space, is called a $primary$ $vector$ $field$. We use $\tau_{+}$ and $\tau_{-}$ to denote the operators which shift the level of descendants, i.e.
$$\tau_{\pm}\Big(\sum_{n,\alpha}f_{n,\alpha}\tau_{n}(\gamma_{\alpha})\Big)
=\sum_{n,\alpha}f_{n,\alpha}\tau_{n\pm1}(\gamma_{\alpha})$$ where $f_{n,\alpha}$ are functions on the big phase space.

Define a $k$-tensor $\gwiig{\cdot\cdot\cdot}$ by
\begin{align}
\gwiig{\mathcal{W}_{1}\mathcal{W}_{2}\cdot\cdot\cdot\mathcal{W}_{k}}
:=\sum_{m_{1},\alpha_{1},...,m_{k},\alpha_{k}}f_{m_{1},\alpha_{1}}^{1}
\cdot\cdot\cdot f_{m_{k},\alpha_{k}}^{k}\frac{\partial^{k}}{\partial t_{m_{1}}^{\alpha_{1}}\cdot\cdot\cdot \partial t_{m_{k}}^{\alpha_{k}}}F_{g}
\end{align}
for vector fields $\mathcal{W}_{i}=\sum_{m,\alpha}f_{m,\alpha}^{i}\frac{\partial}{\partial t_{m}^{\alpha}}$ where $f_{m,\alpha}^{i}$ are functions on the big phase space. This tensor is called $k$-$point$ $(correlation)$ $function$.
For any vector fields $\mathcal{W}_{1}$ and $\mathcal{W}_{2}$ on the big phase space, the $quantum$ $product$ of $\mathcal{W}_{1}$ and $\mathcal{W}_{2}$ is defined by
$$\mathcal{W}_{1}\circ \mathcal{W}_{2}:=\gwii{\mathcal{W}_{1}\mathcal{W}_{2}
\gamma^{\alpha}}\gamma_{\alpha}.$$
This is a commutative and associative product. But it does not have an identity.

Let $$\mathcal{X}:=-\sum_{m,\alpha}(m + d_{\alpha})
\tilde{t}_{m}^{\alpha}\tau_{m}(\gamma_{\alpha})
-\sum_{m,\alpha,\beta}\mathcal{C}_{\alpha}^{\beta}\tilde{t}_{m}^{\alpha}
\tau_{m-1}(\gamma_{\beta})$$
 be the $Euler$ $vector$ $field$ on the big phase space,
where $\tilde{t}_{m}^{\alpha}=t_{m}^{\alpha}-\delta_{m,1}\delta_{\alpha,1}$ and the constants
$d_{\alpha}$ and $\mathcal{C}_{\alpha}^{\beta}$ are uniquely determined by the fact that
the restriction of $\mathcal{X}$  to the small phase space is the $Euler$ $vector$ $field$ of the Frobenius manifold
structure on $\mathcal{H}$.
For example, for quantum cohomology of any compact symplectic manifold $M$, $\mathcal{H}=H^{*}(M,\mathbb{C})$,
$ d_{\alpha}=\frac{1}{2}( \mbox{dimension of }  \gamma_{\alpha})
                - \frac{1}{2}( \mbox{dimension of }  \gamma_{1}) -1, $
 and the matrix $\mathcal{C}=(\mathcal{C}_{\alpha}^{\beta})$ is defined by $c_{1}(M)\cup
 \gamma_{\alpha}=\mathcal{C}_{\alpha}^{\beta}\gamma_{\beta}$ where $c_1(M)$ is the first
 Chern class of $M$.

The quantum multiplication by $\mathcal {X}$ is an endomorphism of the space of primary vector fields on the big phase space. If this endomorphism has distinct eigenvalues at generic points, the big phase space is called $semisimple$. In this case, let $\mathcal{E}_{1},...,\mathcal{E}_{N}$ be the eigenvectors with corresponding eigenvalues $u_{1},...,u_{N}$, i.e.
$$\mathcal {X} \circ \mathcal {E}_{i}=u_{i}\mathcal {E}_{i}$$
for each $i=1,...,N$. $\mathcal{E}_{i}$ is considered as a vector field on the big phase space, and $u_{i}$ is considered as a function on the big phase space. They satisfy the following properties:
$$\mathcal{E}_{i} \circ \mathcal{E}_{j}=\delta_{ij}\mathcal{E}_{i},\quad [\mathcal{E}_{i},\mathcal{E}_{j}]=0,\quad   \mathcal{E}_{i}u_{j}=\delta_{ij}$$
for any $i$ and $j$. Vector fields ${\mathcal{E}_{1},...,\mathcal{E}_{N}}$ are called {\it idempotents} on the big phase space.
When restricted to the small phase space, $\{u_1, \ldots, u_N\}$ gives the
{\it canonical coordinate} system on semisimple Frobenius manifold $\mathcal{H}$ and ${\mathcal{E}_{1},...,\mathcal{E}_{N}}$ coincide with coordinate vector fields of this system (cf. \cite{D}).

Let
$$\mathcal{S}:=-\sum_{m,\alpha}\tilde{t}_{m}^{\alpha}\tau_{m-1}(\gamma_{\alpha})$$
be the $string$ $vector$ $field$ on the big phase space.
For any vector fields $\mathcal{W}$ and $\mathcal{V}$ on the big phase space, define
$$<\mathcal{W},\mathcal{V}>:=\gwii{\mathcal{S}\mathcal{W}\mathcal{V}}.$$
This bilinear form generalizes the pairing $\eta$ on the small phase space. It is non-degenerate only when restricted to the space of primary vector fields.
Define
$$g_{i}:=<\mathcal{E}_{i},\mathcal{E}_{i}>,\quad
h_{i}=\sqrt{g_{i}}$$ Since $<\mathcal{E}_{i},\mathcal{E}_{j}>=0$ if $i\neq j$, functions $g_{1},....,g_{N}$ completely determine the pairing $<\cdot, \,\, \cdot>$ in semisimple case.
When representing expressions involving flat primary fields $\gamma_{\alpha}$ in terms of idempotents, the following simple fact is very useful (cf. \cite{L06}): For any tensor $Q$,
\begin{align}
Q(\gamma_{\alpha},\gamma^{\alpha},\cdot\cdot\cdot)=\sum_{i=1}^{N}
\frac{1}{g_{i}}Q(\mathcal{E}_{i},\mathcal{E}_{i},\cdot\cdot\cdot).
\label{eq:Q}
\end{align}

Define the {\it rotation coefficients on the big phase space} by :
\begin{align}
r_{ij}:=\frac{\mathcal{E}_{j}}{\sqrt{g_{j}}}\sqrt{g_{i}}. \label{r_{ij}}
\end{align}
These functions were introduced in \cite{L06} to study relations among $k$-point functions in semisimple case.
For $i\neq j$, the restriction of $r_{ij}$ to the small phase space coincides with the definition of rotation coefficients for semisimple Frobenius manifolds defined in \cite{D}. This is no longer true when $i=j$.
We will use $\gamma_{ij}$ to denote the rotation coefficients defined in \cite{D}.
In \cite{D}, $\gamma_{ii}$ is defined to be $0$. But ${r_{ii}}$ may not be $0$ according to our definition.
Using the definition of $r_{ij}$ in equation \eqref{r_{ij}}, the
so called {\it string equation} on the small phase space $\mathcal{H}$
has the form
\begin{align}
\sum_{j}r_{ij}h_{j}=0 \label{eq:string}
\end{align} for all $i$.

Now we briefly review basic properties of rotation coefficients which will be used later. Readers are referred to \cite{L06} for more details.
First, $r_{ij}=r_{ji}$, for all $i,j$.
Define
\begin{align*}
v_{ij}:=(u_{j}-u_{i})r_{ij}.
\end{align*}
When taking derivatives of rotation coefficients, the following functions naturally appear:
$$
\theta_{ij}:=\frac{1}{u_{j}-u_{i}}\Big(r_{ij}+\sum_{k}{r_{ik}v_{jk}}
\Big)
$$ for $i\neq j$.
Obviously,
\begin{align}\theta_{ij}+\theta_{ji}=-\sum_{k}r_{ik}r_{jk} \label{eq:theta}\end{align} for any $i\neq j$.
First derivatives of rotation coefficients are given by the formula
\begin{equation}
\mathcal{E}_{k}r_{ij}=r_{ik}r_{jk}+
\begin{cases}
0& \text{if $i\neq j\neq k$},\\
\theta_{ij}& \text{if $k=i\neq j$},\\
\sqrt{\frac{g_{k}}{g_{i}}}\theta_{ik}& \text{if $i=j\neq k$},\\
-2\sum_{l}{r_{il}^{2}}+\sum_{p\neq i}{\sqrt{\frac{g_{p}}{g_{i}}}\theta_{pi}}
+\frac{1}{g_{i}}<\tau_{-}^{2}(\mathcal{S}),\mathcal{E}_{i}>,& \text{if $i=j=k$}.  \label{eq:der of r}
\end{cases}
\end{equation}
Notice that on the small phase space $<\tau_{-}^{2}(\mathcal{S}),\mathcal{E}_{i}>=0$.

Define
$$\Omega_{ij}:
=\frac{1}{u_{j}-u_{i}}\Big(\theta_{ij}-\theta_{ji}
+\sum_{k,l}r_{il}r_{jk}v_{kl}\Big)$$ for $i\neq j$. Then
$$\Omega_{ij}=\Omega_{ji}$$ for all $i\neq j$.
These functions arise naturally in the second order derivatives of rotation coefficients because
$$\mathcal{E}_{j}\theta_{ij}=\Big(r_{jj}-\sqrt{\frac{g_{j}}{g_{i}}}r_{ij}\Big)
\theta_{ij}-\Omega_{ij}$$ for $i\neq j$.
We might consider $\theta_{ij}$ and $\Omega_{ij}$ as functions having poles of order 1 and 2 respectively in terms of $u_{1},...,u_{N}$.

As in \cite{L07}, we will use the following notations for genus-0 and genus-1 $k$-point functions:
$$z_{i_{1},\cdot\cdot\cdot,i_{k}}
:=\gwii{\mathcal{E}_{i_{1}},\cdot\cdot\cdot,\mathcal{E}_{i_{k}}}$$
and $$\phi_{i_{1},\cdot\cdot\cdot,i_{k}}:=
\gwiione{\mathcal{E}_{i_{1}},\cdot\cdot\cdot,\mathcal{E}_{i_{k}}}.$$
In \cite{L06}, it was proved that genus-0 4-point functions satisfy the following properties: For $i\neq j$,
\begin{align}
(i) &\quad z_{iiii}=-g_{i}r_{ii},
\nonumber\\
(ii)&\quad  z_{jiii}=-z_{jjii}=-\sqrt{g_{i}g_{j}}r_{ij},
\nonumber\\
(iii) &\quad z_{ijkl}=0, \quad \mbox{otherwise}. \label{eq:g=0 4-pt}
\end{align}
It was also proved that the genus-1 1-point functions are given by
\begin{align} 24\phi_{i}=-12\sum_{j}{r_{ij}v_{ij}}
-\sum_{j}{\frac{h_{i}}{h_{j}}}r_{ij}
\label{eq:g=1 1-pt}
\end{align} for all $i$.
Higher point genus-0 and genus-1 functions can be computed recursively using the following formula:
\begin{align}
&\gwiig{\mathcal{E}_{i_{1}}\mathcal{E}_{i_{2}}\cdot\cdot\cdot
\mathcal{E}_{i_{k+1}}}\nonumber
\\&=\mathcal{E}_{i_{k+1}}\gwiig{\mathcal{E}_{i_{1}}\cdot\cdot\cdot
\mathcal{E}_{i_{k}}}-\Big(\sum_{j=1}^{k}r_{i_{j},i_{k+1}}
\sqrt{\frac{g_{i_{k+1}}}{g_{i_{j}}}} \Big)\gwiig{\mathcal{E}_{i_{1}}\cdot\cdot\cdot\mathcal{E}_{i_{k}}}
\nonumber
\\&-\sum_{j=1}^{k}r_{i_{j},i_{k+1}}
\sqrt{\frac{g_{i_{j}}}{g_{i_{k+1}}}}
\gwiig{\mathcal{E}_{i_{1}}\cdot\cdot\cdot\widehat{\mathcal{E}_{i_{j}}}
\cdot\cdot\cdot\mathcal{E}_{i_{k}}\mathcal{E}_{i_{k+1}}}
\nonumber
\\&+\sum_{j=1}^{k}\delta_{i_{k+1},i_{j}}\sum_{p}r_{p,i_{k+1}}
\sqrt{\frac{g_{i_{k+1}}}{g_{p}}}\gwiig{\mathcal{E}_{i_{1}}\cdot\cdot\cdot
\widehat{\mathcal{E}_{i_{j}}}
\cdot\cdot\cdot\mathcal{E}_{i_{k}}\mathcal{E}_{p}}
\label{eq:ditui}
\end{align}
for all $g \geq 0$.
For example, the genus-1 2-point functions are given by
\begin{align}
24\phi_{ij}=& 12r_{ij}^2 + \sum\limits_l {({r_{il}}{r_{jl}}{{{h_i}{h_j}} \over {h_l^2}} - {r_{ij}}{r_{il}}{{{h_j}} \over {{h_l}}} - {r_{ij}}{r_{jl}}{{{h_i}} \over {{h_l}}})}-\{{\theta _{ij}}{{{h_j}} \over {{h_i}}} + {\theta _{ji}}{{{h_i}} \over {{h_j}}}\}
\nonumber\\&-24r_{ij}\{\frac{h_{j}}{h_{i}}\phi_{i}+\frac{h_{i}}{h_{j}}\phi_{j}\}
\label{eq:g=1 2-ptij}
\end{align}
for $i\neq j$, and
\begin{align}
24\phi_{ii}=&12r_{ii}^{2}+\sum_{j}\{r_{ij}^{2}(-10+\frac{h_{i}^{2}}{h_{j}^{2}})
-2r_{ii}r_{ij}\frac{h_{i}}{h_{j}}\}-\sum_{j\neq i}\{\theta_{ij}\frac{h_{i}}{h_{j}}+\theta_{ji}\frac{h_{j}}{h_{i}}\}
\nonumber
\\&-48r_{ii}\phi_{i}+24\sum_{j}r_{ij}\frac{h_{i}}{h_{j}}\phi_{j}-\frac{1}{h_{i}^{2}}
<\tau_{-}^{2}(S),\mathcal{E}_{i}>
\label{eq:g=1 2-ptii}
\end{align} for all $i$.
Some genus-1 3-point functions can be found in Appendix~\ref{sec:3ptg1}.

For any vector field $\mathcal{W}$, define
$$T(\mathcal{W}):=\tau_{+}(\mathcal{W})-\gwii{\mathcal{W}\gamma^{\alpha}}
\gamma_{\alpha}.$$ The operator $T$ was introduced in \cite{L02} in order to simplify topological recursion relations for Gromov-Witten invariants. The following properties of operator $T$ were proved in \cite{L06}:
\begin{align}
(i)&\quad  T(\mathcal{E}_{j}) \, u_{i}=0,  \nonumber \\
(ii)&\quad T(\mathcal{E}_{j}) \, g_{i}=-2\delta_{ij}g_{i}, \nonumber \\
(iii)&\quad  T(\mathcal{E}_{k}) \,r_{ij}=\delta_{ij}r_{ik}\sqrt{\frac{g_{k}}{g_{i}}}, \nonumber \\
(iv)&\quad T(\mathcal{E}_{k}) \, \theta _{ij} =- {r_{ij}}{r_{ik}}\sqrt {\frac{{{g_k}}}{{{g_i}}}} \quad\mbox{for $i \ne j$},
        \nonumber \\
(v)& \quad T(\mathcal{E}_{k}) \, \Omega _{ij}
    = \theta _{ji} {r_{ik}}\sqrt {\frac{{{g_k}}}{{{g_i}}}}
    + \theta _{ij}  {r_{jk}}\sqrt {\frac{{{g_k}}}{{{g_j}}}}  \quad\mbox{for $i \ne j$},  \nonumber \\
(vi)& \quad  T(\mathcal{E}_{k}) < \tau _ - ^2(S),\mathcal{E}_{i} >  =  - \delta_{ik} < \tau _ - ^2(S),\mathcal{E}_{i} >,  \nonumber \\
(vii)& \quad T(\mathcal{E}_{k}) < \tau _ - ^3(S),\mathcal{E}_{i} >  =  - \delta_{ik} < \tau _ - ^3(S),\mathcal{E}_{i} >.
\label{eqn:Tei}
\end{align}

\section{Proof of the main theorem}
\label{sec:pfmain}
\allowdisplaybreaks

 Recall that for a semisimple Frobenius manifold $\mathcal{H}$, the genus-2 G-function  $G^{(2)}(u,u_{x},u_{xx})$ defined in \cite{DLZ} has the following form:
                  \begin{align}
                   G^{(2)}(u,u_{x},u_{xx})
                   =&\sum_{i}G_{i}^{(2)}(u,u_{x})u_{xx}^{i}
                   +\sum_{i\neq j}G_{ij}^{(2)}(u)\frac{(u_{x}^{j})^{3}}{u_{x}^{i}}
                   \nonumber\\&+\frac{1}{2}\sum_{i,j}P_{ij}^{(2)}(u)u_{x}^{i}u_{x}^{j}
                   +\sum_{i}Q_{i}^{(2)}(u)(u_{x}^{i})^{2}
                              \end{align}
 where $G_{i}^{(2)}$, $G_{ij}^{(2)}$, $P_{ij}^{(2)}$, and $Q_{i}^{(2)}$ are functions on the jet space of $\mathcal{H}$. Precise
 definitions for these functions can be found in Appendix~\ref{sec:G2}.
In this formula, $u_{i}=u^{i}$ for $i=1,...,N$ are the canonical coordinates on the small phase space
$\mathcal{H}$, i.e.
the restriction of functions $u_{i}$ defined in Section~\ref{sec:prel} to $\mathcal{H}$.
 Let  $w^{\alpha}$ for $\alpha=1,...,N$ be the flat coordinates on $\mathcal{H}$ with respect to
 the basis $\{\gamma_{1},\cdot\cdot\cdot,\gamma_{N}\}$, i.e. $w^{\alpha}$ is the restriction of
 $t_0^{\alpha}$ to the small phase space $\mathcal{H}$. Under the transformation
 \begin{equation} \label{eqn:s2b}
 w^{\alpha}=\gwii{\gamma_{1}\gamma^{\alpha}},
 \end{equation}
 every function on the small phase space can also be viewed as a function on the big phase space.
 In particular,  $u_{i,x}=u^{i}_{x}$ means the first order derivative of $u_{i}$ with respect to $x:=t_{0}^{1}$
 after transformation \eqref{eqn:s2b}.
 The functions $u_{i},u_{i,x},u_{i,xx},...$ for $i=1,...,N$ form a coordinate system for the jet space of the semisimple Frobenius manifold $\mathcal{H}$ (cf: \cite{DZ}).

To prove that the genus-2 G-function is 0, it suffices to show that functions $G_{i}^{(2)}$ and
$Q_{i}^{(2)} + \frac{1}{2} P_{ii}^{(2)}$
for all $i$, $G_{ij}^{(2)}$ and
$P_{ij}^{(2)} + P_{ji}^{(2)} $ for $i \neq j$,
are equal to 0.
In this section, we prove that all these functions  vanish under conditions $(C1)$, $(C2)$ and $(C3)$.
The main idea of the proof is to first represent all these functions in terms of rotation coefficients $r_{ij}$ and
functions $h_{i}$, $u_i$, $\theta_{ij}$, and $\Omega_{ij}$ defined in Section~\ref{sec:prel}. We then use conditions
$(C1)$, $(C2)$ and $(C3)$ to
get rid of $u_i$, $\theta_{ij}$, and $\Omega_{ij}$ and obtain expressions for the above four types of functions only
involving $r_{ij}$ and $h_{i}$. Then straightforward calculations show that all these functions vanish.

The following lemma will be very useful in the proof.
\begin{lem} Under the condition $(C2)$, on the small phase space the following identities hold:
\begin{align}
&\sum_{j}{r_{ij}v_{ij}} \,\, = \,\,
-\frac{1}{12}\sum_{j}r_{ij}\frac{h_{i}}{h_{j}}
\label{eq:g=1 1-pt small}
\end{align}
for all $i$,
and
\begin{align}
&{\theta _{ij}}{{{h_j}} \over {{h_i}}} + {\theta _{ji}}{{{h_i}} \over {{h_j}}} \,\, = \,\, 12r_{ij}^2 + \sum\limits_l {({r_{il}}{r_{jl}}{{{h_i}{h_j}} \over {h_l^2}} - {r_{ij}}{r_{il}}{{{h_j}} \over {{h_l}}} - {r_{ij}}{r_{jl}}{{{h_i}} \over {{h_l}}})}\label{eq:g=1 2-pt small}
  \end{align}
  for all $i\neq j$.
\end{lem}
{\bf Proof}:
First note that condition $(C2)$ is equivalent to the vanishing of the restriction of genus-1 correlation functions
to the small phase space, i.e.
\begin{align}
\phi_{i_{1}\cdot\cdot\cdot i_{k}}\mid_{\mathcal{H}}=0 \label{eq:g=1 k-pt}
\end{align}
for any
$\{i_{1},\cdot\cdot\cdot,i_{k}\}$, with $k\geq 1$.

Equation $(\ref{eq:g=1 1-pt small})$ follows from equation (\ref{eq:g=1 1-pt}) and $\phi_{i}\mid_{\mathcal{H}}=0$.
Equation (\ref{eq:g=1 2-pt small}) follows from equation (\ref{eq:g=1 2-ptij}),
$\phi_{i}\mid_{\mathcal{H}}=0$, and
$\phi_{ij}\mid_{\mathcal{H}}=0$ for  $i\neq j$.
$\Box$

In the definition of the genus-2 G-function (see Appendix~\ref{sec:G2}),
the following functions were introduced to simplify expressions:
\begin{align*}
H_{i}:=\frac{1}{2}\sum_{j\neq i}u_{ij}\gamma_{ij}^{2}
\end{align*}
where $u_{ij}:=u_{i}-u_{j}$ and $\gamma_{ij}$ are the rotation coefficients defined in \cite{D}.
An immediate application of equation (\ref{eq:g=1 1-pt small})
is to get rid of functions $u_{ij}$ in $H_i$, i.e.
\begin{align}
H_{i}=\frac{1}{24}\sum_{k}r_{ik}\frac{h_{i}}{h_{k}}
\label{eq:H}
\end{align}
when condition (C2) is satisfied.

\subsection{ Vanishing of $G_{ij}^{(2)}$}

In this subsection, we prove $G_{ij}^{(2)}=0$ for all $i\neq j$ under condition (C2).
In the definition of $G_{ij}^{(2)}$ given in Appendix \ref{sec:G2},
$\partial_{i} := \frac{\partial}{\partial u_{i}}$ is equal to taking derivative along
direction $\mathcal{E}_{i}$ on the small phase space. Therefore we can compute $\partial_i h_j$ and
$\partial_i \gamma_{jk}$   using formulas \eqref{r_{ij}} and  \eqref{eq:der of r}.
We can also replace $H_i$ by the  right hand side of equation~\eqref{eq:H} . We then obtain the
following expression for $G_{ij}^{(2)}$:
  \begin{align*}
 G_{ij}^{(2)}
 =-\frac{r_{ij}}{5760h_{i}h_{j}}\big\{{\theta _{ij}}{{{h_j}} \over {{h_i}}} + {\theta _{ji}}{{{h_i}} \over {{h_j}}}- 12r_{ij}^2 - \sum\limits_l {({r_{il}}{r_{jl}}{{{h_i}{h_j}} \over {h_l^2}} - {r_{ij}}{r_{il}}{{{h_j}} \over {{h_l}}} - {r_{ij}}{r_{jl}}{{{h_i}} \over {{h_l}}})}\big\}.
 \end{align*}
The right hand side of this equation is understood as a function obtained by applying transformation \eqref{eqn:s2b} to
a function on the small phase space. Therefore by equation (\ref{eq:g=1 2-pt small}), we have $ G_{ij}^{(2)}=0$.
Most formulas in the proof of Theorem~\ref{thm:main}  will be understood in a similar way. To make the paper concise, we will not repeat this argument later.

This simple proof of $ G_{ij}^{(2)}=0$ is rather typical in the rest part of the proof for Theorem~\ref{thm:main}.
To prove the vanishing of other parts of the genus-2 G-function, we will use similar ideas although the computations
will be much more involved.


\subsection{ Vanishing of $G_{i}^{(2)}$}
\label{sec:Gi}

 In this subsection, we show $G_{i}^{(2)}=0$ for all $i$ under conditions (C1) and (C2).
 We write
 \[ G_{i}^{(2)} = G_{i, 1}^{(2)} + G_{i, 2}^{(2)} \]
 where $G_{i, 1}^{(2)}$ explicitly involves jet coordinates $u_{k, x}$ but $G_{i, 2}^{(2)}$ does not.
 Precise definitions for $G_{i, 1}^{(2)}$ and $G_{i, 2}^{(2)}$ are given in Appendix \ref{sec:G2}.
 We will compute these two functions separately.

 First note that for any function $f$ obtained from applying transformation~\eqref{eqn:s2b} to
 a function on the small phase space,
 \[ \partial_x f = \sum_{k} u_{k , x} \partial_{k} f \]
 by the chain rule.
Using formulas \eqref{r_{ij}} and  \eqref{eq:der of r}, we can then compute
all partial derivatives involved in the  definitions of  $G_{i, 1}^{(2)}$ and $G_{i, 2}^{(2)}$.
 We can also compute $H_i$ using
 equation \eqref{eq:H}. After straightforward calculations, we obtain
 the following formulas:
\begin{align*}
G_{i,1}^{(2)}=&\sum_{k\neq i}\frac{1}{1920h_{i}^{2}}\big\{{\theta _{ik}}{{{h_k}} \over {{h_i}}} + {\theta _{ki}}{{{h_i}} \over {{h_k}}}- 12r_{ik}^2 - \sum\limits_l {({r_{il}}{r_{kl}}{{{h_i}{h_k}} \over {h_l^2}} - {r_{ik}}{r_{il}}{{{h_k}} \over {{h_l}}} - {r_{ik}}{r_{kl}}{{{h_i}} \over {{h_l}}})}\big\}\frac{\partial_{x}u_{k}}{\partial_{x}u_{i}}
\\&+\sum_{k\neq
i}\frac{\theta_{ik}}{5760h_{i}h_{k}}-\frac{r_{ii}^2}{1440h_{i}^{2}}-\sum_{k}
\frac{r_{ik}^{2}}{1920h_{k}^{2}}+\sum_{k}\frac{r_{ii}r_{ik}}{1440h_{i}h_{k}}
\end{align*}
 and
\begin{align}
G_{i,2}^{(2)}
=&\sum_{k\neq i}\frac{\theta_{ik}}{2880h_{i}h_{k}}-\sum_{k\neq i}\frac{\theta_{ki}h_{k}}{384h_{i}^{3}}
    +\sum_{k\neq i}\frac{7\theta_{ki}}{2880h_{i}h_{k}}-\sum_{j}\frac{17r_{ii}r_{ij}}{2880h_{i}h_{j}}
    +\sum_{k}\frac{r_{ik}r_{kk}h_{i}}{1440h_{k}^{3}} \nonumber
\\&-\sum_{k}\frac{19r_{ik}^{2}}{720h_{i}^{2}}+\sum_{k}\frac{r_{ik}^{2}}{1440h_{k}^{2}}
    +\sum_{j,k}\frac{r_{ik}r_{ij}}{2880h_{j}h_{k}}
    +\sum_{j,k}\frac{7r_{ik}r_{jk}}{2880h_{i}h_{j}}-\sum_{k,l}\frac{h_{i}r_{il}r_{kl}}{2880h_{k}h_{l}^{{2}}} \nonumber
\\&+\frac{23r_{ii}^{2}}{720h_{i}^{2}}. \label{eqn:Gi2}
\end{align}
By equation (\ref{eq:g=1 2-pt small}), the coefficient for
$\frac{\partial_{x}u_{k}}{\partial_{x}u_{i}}$ in $G_{i,1}^{(2)}$ becomes $0$ under condition (C2).
So we have
\begin{align}
G_{i,1}^{(2)}=&\sum_{k\neq
i}\frac{\theta_{ik}}{5760h_{i}h_{k}}-\frac{r_{ii}^2}{1440h_{i}^{2}}-\sum_{k}
\frac{r_{ik}^{2}}{1920h_{k}^{2}}+\sum_{k}\frac{r_{ii}r_{ik}}{1440h_{i}h_{k}}.
\label{eqn:Gi1}
\end{align}

Next we will get rid of first order poles $\theta_{ik}$ and $\theta_{ki}$ in equations \eqref{eqn:Gi1} and \eqref{eqn:Gi2}
 using functions only depending on $\{ r_{ij} \mid i, j = 1, \ldots , N \}$ and $\{ h_{i} \mid i=1, \ldots, N  \}$.
 Here we also need condition (C1).
 A key ingredient in this process is the following lemma:
 \begin{lem} \label{lem:Gi}
Under conditions $(C1)$ and $(C2)$, on the small phase space we have
   \begin{align}
&2\sum\limits_{\mathop i\limits_{i \ne k} } \frac{\theta_{ki}}{h_i} =
\sum\limits_i {(7{{h_k} \over {h_i^2}}r_{ik}^2 - 6 \frac{r_{ik}^2}{h_k}
- 2{{h_k^2} \over {h_i^3}}{r_{ii}}{r_{ik}})}
+ \sum\limits_{i,j} {({{h_k^2} \over {{h_i}h_j^2}}{r_{ij}}{r_{jk}}
- {{h_k} \over {{h_i}{h_j}}}{r_{ik}}{r_{jk}})}
\label{eq:sumtheta}
   \end{align}
for all $k$.
\end{lem}

\noindent{\bf Proof: }
Condition (C1) is equivalent to say that the restriction of
$\gwii{ {\gamma _\alpha }{\gamma ^\alpha }{\gamma _\beta }{\gamma ^\beta } }$
to the small phase space is constant.
   By equation (\ref{eq:Q}), we have
\begin{align*}
  { \gwii{ {\gamma _\alpha }{\gamma ^\alpha }{\gamma _\beta }{\gamma ^\beta }{ }}}
  &=\sum_{i,j=1}^{N}\frac{1}{g_{i}g_{j}}
   \gwii{\mathcal{E}_{i}\mathcal{E}_{i}\mathcal{E}_{j}\mathcal{E}_{j}}
  =\sum_{i,j=1}^{N}\frac{1}{g_{i}g_{j}}z_{iijj}.
\end{align*}
So equation \eqref{eq:g=0 4-pt} implies
\begin{equation}
\gwii{ {\gamma _\alpha }{\gamma ^\alpha }{\gamma _\beta }{\gamma ^\beta } }
 \, = \, \sum_{i,j}\frac{r_{ij}}{h_{i}h_{j}}-2\sum_{i}\frac{r_{ii}}{h_{i}^{2}}.
   \label{eq:O}
\end{equation}
Therefore if condition (C1) is satisfied, on the small phase space, we have
\[ \sum_{i,j}\frac{r_{ij}}{h_{i}h_{j}}-2\sum_{i}\frac{r_{ii}}{h_{i}^{2}} \, = \, \mbox{constant}. \]
Taking derivative of this equation along the direction of  $\mathcal{E}_{k}$, we have
\begin{align}
\sum\limits_{i \ne k} ( {{{h_k}} \over {h_i^3}} + {{{h_i}} \over {h_k^3}}){\theta _{ik}} - 2\sum\limits_{i \ne k} {{{{\theta _{ki}}} \over {{h_k}{h_i}}}}
= &\sum\limits_{i,j} {{{{r_{ik}}{r_{jk}}} \over {{h_i}{h_j}}}}  - 2\sum\limits_i {{{r_{ik}^2} \over {h_i^2}}}  - \sum\limits_{i,j} {{{{r_{ij}}{h_k}({r_{ik}}{h_j} + {r_{jk}}{h_i})} \over {h_i^2h_j^2}}} \nonumber \\
&+ 4\sum\limits_i {{{{h_k}{r_{ik}}{r_{ii}}} \over {h_i^3}}}  + 2\sum\limits_i {{{r_{ik}^2} \over {h_k^2}}}.
\label{eqn:derO1-O2}
\end{align}
We then convert the first summation on the left hand side to an expression similar to
the second summation on the left hand side together with some function only depending on $r_{ij}$ and $h_i$.
To see how this can be done, we write
\begin{align*}
&\sum_{i\neq k}\frac{h_{k}}{h_{i}^{3}}\theta_{ik}
=\sum_{i\neq k}\{\frac{h_{k}}{h_{i}}\theta_{ik}\}\frac{1}{h_{i}^{2}}
\end{align*}
and use equation (\ref{eq:g=1 2-pt small}) to replace $\frac{h_{k}}{h_{i}}\theta_{ik}$
by $-\frac{h_{i}}{h_{k}}\theta_{ki}$ plus a function only depending on $r_{ij}$ and $h_i$.
In this way, we can see that
\begin{align}
&\sum_{i\neq k}\frac{h_{k}}{h_{i}^{3}}\theta_{ik}
=-\sum_{i\neq k}\{\frac{1}{h_{i}h_{k}}\theta_{ki}\}+\{\mbox{some function of $r$ and $h$}\}.
\label{eqn:theta/3}
\end{align}
Here ``some function of $r$ and $h$" means a function which only depends on $h_i$ for all $i$ and
 $r_{ij}$ for all $i$ and $j$.  For simplicity, we omit the precise form of such functions.
Note that we can use equation (\ref{eq:theta}) to switch indices in $\theta_{ik}$. So by a similar calculation using
equation (\ref{eq:g=1 2-pt small}), we also obtain
\begin{align*}
&\sum_{i\neq k}\frac{h_{i}}{h_{k}^{3}}\theta_{ik}
= -\sum_{i\neq k}\{\frac{1}{h_{i}h_{k}}\theta_{ki}\}+\{\mbox{some function of $r$ and $h$}\}.
\end{align*}
After plugging these terms into equation~\eqref{eqn:derO1-O2} and multiplying the resulting expression by
$- \frac{1}{2} h_k$, we obtain equation (\ref{eq:sumtheta}). $\Box$

\begin{rem} \label{rem:Gi}
Combining equations $(\ref{eq:theta})$ and $(\ref{eq:sumtheta})$,
we can  express $\sum\limits_{\mathop i\limits_{i \ne k} } \frac{\theta_{ik}}{h_i}$
in terms of functions $\{ r_{ij} \}$ and $\{h_i \}$.
By equations $(\ref{eq:theta})$, $(\ref{eq:g=1 2-pt small})$ and $(\ref{eq:sumtheta})$,
 we can also express $\sum\limits_{\mathop i\limits_{i \ne k} } h_i \theta _{ki}$
  and $\sum\limits_{\mathop i\limits_{i \ne k} } h_i \theta _{ik}$ in terms of functions
  $\{ r_{ij} \}$ and $\{ h_i \}$. More generally,
  we can express both $\sum\limits_{\mathop i\limits_{i \ne k} } h_i^p \theta _{ik}$
  and $\sum\limits_{\mathop i\limits_{i \ne k} } h_i^p \theta _{ki}$ in terms of functions
  $\{ r_{ij} \}$ and $\{ h_i \}$ if $p$ is odd. We also note that repeatedly applying equations
  \eqref{eq:theta} and \eqref{eq:g=1 2-pt small} shows that
 \[ \theta_{ij} h_i^p h_j^q \,\, = \,\, \theta_{ij} h_i^{p+2k} h_j^{q-2k} + \{\mbox{some function of $r$ and $h$}\}\]
 for any integers $p$, $q$, and $k$. This observation is very useful in the process of getting
 rid of functions $\theta_{ij}$ in complicated expressions.

\end{rem}

Now we come back to the computation of $G_i^{(2)}=G_{i, 1}^{(2)} + G_{i, 2}^{(2)}$. After replacing
all terms containing $\theta_{ik}$ and $\theta_{ki}$ in equations \eqref{eqn:Gi1} and \eqref{eqn:Gi2}
using the method indicated in Lemma~\ref{lem:Gi} and Remark~\ref{rem:Gi}, we obtain a formula
for  $G_{i}^{(2)}$ which only involves functions $\{r_{jk} \mid j, k=1, \ldots, N\}$ and $\{h_j \mid j= 1, \ldots, N\}$.
Straightforward computations
then show  $G_{i}^{(2)}=0$.

\subsection{ Vanishing of $P_{ij}^{(2)}+P_{ji}^{(2)}$}
\label{sec:P=0}

In this subsection, we prove $P_{ij}^{(2)}+P_{ji}^{(2)}=0$ for all $i \neq j$ under conditions (C1) and (C2).
We can use formulas \eqref{r_{ij}} and  \eqref{eq:der of r} to compute
all partial derivatives involved in the definition of $P_{ij}^{(2)}$ given in Appendix~\ref{sec:G2}.
We can also compute functions $H_i$ using equation \eqref{eq:H}. We then obtain a formula for
$P_{ij}^{(2)}$ only depending on $\{r_{kl} \mid k, l=1, \ldots, N\}$, $\{h_k \mid k=1, \ldots, N\}$,
and first order poles $\{ \theta_{kl} \mid k, l=1, \ldots, N\}$. The precise form of this formula
is given by equation~\eqref{eqn:Pijc} in  Appendix~\ref{sec:PQ}.
Note that we have used equation~\eqref{eq:theta} to get rid of some
terms containing $\theta_{ik}$ when deriving this formula.
As in subsection~\ref{sec:Gi}, the main idea for proving $P_{ij}^{(2)}+P_{ji}^{(2)}=0$ is to replace
first order poles $\theta_{kl}$ in $P_{ij}^{(2)}+P_{ji}^{(2)}$ by functions only depending on $\{r_{kl} \}$
and $\{h_k\}$. For simplicity of notations, $\{r_{kl} \}$ means the set of all functions
$r_{kl}$ where $k$ and $l$ run over their entire  range. The notation $\{h_k\}$ is interpreted
in a similar way.

For convenience, we give all terms in $P_{ij}^{(2)}$ which contain first order poles below:
 \begin{align}
P_{ij}^{(2)}=
&\frac{1}{1440} \left( -3 \sum_{k\neq j}\frac{h_{i}h_{j}r_{ik}\theta_{kj}}{h_{k}^{4}}
    -\frac{41r_{ij}\theta_{ij}}{h_{i}^{2}}
+\frac{6r_{ii}h_{j}\theta_{ij}}{h_{i}^{3}}
-3 \sum_{k}\frac{r_{ik}h_{j}\theta_{ij}}{h_{i}^{2}h_{k}} \right.
\nonumber
\\& \left. \hspace{60pt} -4 \sum_{k\neq i}\frac{r_{ij}h_{j}h_{k}\theta_{ik}}{h_{i}^{4}}
- 5 \sum_{k\neq i}\frac{r_{ij}h_{j}\theta_{ki}}{h_{i}^{2}h_{k}} \right)
\nonumber
\\&+\{\mbox{some function of $r$ and $h$}\}.
\label{eqn:Pijp}
\end{align}
Using Remark~\ref{rem:Gi},
the two terms containing $\theta_{ik}$ and $\theta_{ki}$
 in the second line of equation~\eqref{eqn:Pijp} can be expressed as functions of $\{ r_{kl} \}$ and $\{ h_k \}$.
Below we explain how to deal with the first 4 terms in equation~\eqref{eqn:Pijp}. For simplicity, we will
omit precise forms of functions only depending on $\{ r_{kl} \}$ and $\{ h_k \}$.

 By equation (\ref{eq:g=1 2-pt small}), we have
         \begin{align}
\sum_{k\neq j}\frac{h_{i}h_{j}r_{ik}\theta_{kj}}{h_{k}^{4}}
= - \frac{h_{i}}{h_{j}}\sum_{k\neq j}\frac{r_{ik}\theta_{jk}}{h_{k}^{2}}
+\{\mbox{some  function  of  $r$  and  $h$}\}. \label{P}
\end{align}
Next we will use the following identities
\begin{lem} \label{lem:thetaijk}
Under condition $(C2)$, on the small phase space we have
\begin{align}
& {\theta _{ik}}{r_{jk}}\frac{{{h_j}}}{{{h_i}}} + {\theta _{ki}}{r_{ij}}\frac{{{h_j}}}{{{h_k}}}
 + {\theta _{jk}}{r_{ik}}\frac{{{h_i}}}{{{h_j}}} + {\theta _{kj}}{r_{ij}}\frac{{{h_i}}}{{{h_k}}}
 + {\theta _{ij}}r_{jk}\frac{{{h_k}{}}}{{{h_i}}} + {\theta _{ji}}{r_{ik}}\frac{{{h_k}}}{{{h_j}}}
                            \nonumber
 \\=&\{\mbox{some  function  of  $r$  and  $h$}\}
 \label{g=1 3-pt jian}
\end{align}
  for  $i\neq j\neq k$. For precise formula, see equation (\ref{g=1 3-ptijkjian}) in Appendix~\ref{sec:PQ}.
\end{lem}
{\bf Proof:}
Condition (C2) implies that genus-1 correlation functions on the small phase space $\mathcal{H}$ are $0$.
For distinct $i$, $j$, $k$, we can compute genus-1 3-point function $\phi_{ijk}$ using
 formula $(\ref{eq:ditui})$ together with  the formulas of genus-1 1-point functions $(\ref{eq:g=1 1-pt})$
and genus-1 2-point functions $(\ref{eq:g=1 2-ptij})$. The precise formula of $\phi_{ijk}$
is given as equation (\ref{g=1 3-ptijk }) in Appendix~\ref{sec:3ptg1}.
The lemma just follows from
$\phi_{i}|_{\mathcal{H}}=0$, $\phi_{ij}|_{\mathcal{H}}=0$, $\phi_{ijk}|_{\mathcal{H}}=0$
for any $i\neq j\neq k$. $\Box$

Multiplying both sides of equation (\ref{g=1 3-pt jian}) by $\frac{1}{h_{k}^{2}}$
and taking sum over $k$ for $k\notin \{i,j\}$,  we  obtain
   \begin{align}
   &\frac{h_{i}}{h_{j}}\sum_{k\neq j}\frac{r_{ik}\theta_{jk}}{h_{k}^{2}}+\frac{h_{j}}{h_{i}}\sum_{k\neq j}\frac{r_{jk}\theta_{ik}}{h_{k}^{2}}
  \nonumber\\=&  - \sum_{k}{\theta _{ij}}r_{jk}\frac{{{1}}}{{{h_i}{h_k}}}
 - \sum_{k}{\theta _{ji}}{r_{ik}}\frac{{1}}{{{h_j}{h_k}}}
+2{\theta _{ij}}{r_{jj}}\frac{{1}}{{{h_i}{h_j}}}
  +2{\theta _{ji}}{r_{ii}}\frac{{1}}{{{h_i}{h_j}}}
   \nonumber\\&+\mbox{\{some function of $r$ and $h$\}}. \label{P1}
   \end{align}
In deriving this formula, we have used equation \eqref{eq:g=1 2-pt small} and the fact that both
$\sum_{k \neq i} \frac{\theta_{ki}}{h_k^3}$ and
$\sum_{k \neq j} \frac{\theta_{kj}}{h_k^3}$
can be expressed as functions of $\{r_{kl} \}$ and $\{h_k\}$, which in turn follow from
equations \eqref{eqn:theta/3} and \eqref{eq:sumtheta}. By equations (\ref{P}) and (\ref{P1}), we get
  \begin{align*}
&P_{ij}^{(2)}+P_{ji}^{(2)}
\\=& \frac{1}{1440} \left( - 3\sum_{k}r_{jk}{\theta _{ij}}\frac{{{1}}}{{{h_i}{h_k}}}
- 3 \sum_{k}r_{jk}\theta_{ji}\frac{h_{i}}{h_{j}^{2}h_{k}}
- 3 \sum_{k}{r_{ik}}{\theta _{ji}}\frac{{1}}{{{h_j}{h_k}}}
-3 \sum_{k}r_{ik}\theta_{ij}\frac{h_{j}}{h_{i}^{2}h_{k}}
 \right.
\\&  \left. \hspace{40pt} - \frac{41 r_{ij}}{h_{i}h_{j}}
(\theta_{ij}\frac{h_{j}}{h_{i}}+\theta_{ji}\frac{h_{i}}{h_{j}})
+6r_{ii}\theta_{ij}\frac{h_{j}}{h_{i}^{3}}
 +6 {r_{ii}}{\theta _{ji}}\frac{{1}}{{{h_i}{h_j}}}
 +6r_{jj}\theta_{ji}\frac{h_{i}}{h_{j}^{3}}
+ 6{r_{jj}}{\theta _{ij}}\frac{{1}}{{{h_i}{h_j}}}
 \right)
\\&
+\mbox{\{some function of $r$ and $h$\}}.
    \end{align*}
Using equation \eqref{eq:g=1 2-pt small}, it is easy to see that the right hand side of this equation
can be expressed as a function only depending on $\{r_{kl}\}$ and $\{h_k\}$.
 A straightforward calculation
shows that this function is equal to $0$. This proves that
$P_{ij}^{(2)}+P_{ji}^{(2)}=0$ for $i \neq j$.

\subsection{ Vanishing of $\frac{1}{2}P_{ii}^{(2)}+Q_{i}^{(2)}$}
\label{sec:P+Q=0}

In this subsection, we prove $\frac{1}{2}P_{ii}^{(2)}+Q_{i}^{(2)}=0$ for all $i$
under conditions (C1), (C2), and (C3). In the definition of $Q_{i}^{(2)}$  given in
Appendix~\ref{sec:G2}, there are 6 terms having $u_{ik} = u_i - u_k$ in denominator.
After calculating derivatives of rotation coefficients by equation \eqref{eq:der of r} and
replacing $H_k$ by equation \eqref{eq:H}, these terms become
\begin{align}
& \frac{1}{1152} \left(  \sum\limits_{\mathop {k,l}\limits_{k \ne i,l \ne i} }
                {\frac{{{h_k}{v_{lk}}{\theta _{il}}}}{{{u_{ik}}h_i^3}}}
    - 2 \sum\limits_{k \ne i} {\frac{{{\theta _{ik}}{h_k}}}{{{u_{ik}}h_i^3}}}
    + \frac{1}{12}\sum\limits_{\mathop {k,l}\limits_{k \ne i} }
                {\frac{{{r_{ik}}{r_{il}}{h_k}}}{{{u_{ik}}h_i^2{h_l}}}}
    - \sum\limits_{k \ne i} {\frac{{{r_{ii}}{r_{ik}}{h_k}}}{{{u_{ik}}h_i^3}}} \right.
    \nonumber
\\& \hspace{28pt} - 2 \sum\limits_{k \ne i} {\frac{{{\theta _{ki}}}}{{{u_{ik}}{h_i}{h_k}}}}
    + \frac{1}{12}\sum\limits_{\mathop {k,l}\limits_{k \ne i} }
                    {\frac{{{r_{ik}}{r_{kl}}}}{{{u_{ik}}{h_i}{h_l}}}}
    - 2 \sum\limits_{k \ne i} {\frac{{{r_{kk}}{r_{ik}}}}{{{u_{ik}}{h_i}{h_k}}}}
    +  \sum\limits_{\mathop {k,l}\limits_{l \ne i} }
                    {\frac{{{r_{ik}}{r_{ll}}{v_{kl}}}}{{{u_{il}}{h_i}{h_l}}}}
                    \nonumber
\\& \hspace{28pt} \left. + \sum\limits_{\mathop {k,l}\limits_{k \ne l,l \ne i} }
                    {\frac{{{u_{lk}}{r_{ik}}{\theta _{lk}}}}{{{u_{il}}{h_i}{h_l}}}}
    + \sum\limits_{l \ne i} {\frac{{{r_{ii}}{\theta _{li}}}}{{{h_i}{h_l}}}}
    +  \sum\limits_{k \ne i} {\frac{{{r_{ii}}{h_k}{\theta _{ik}}}}{{h_i^3}}}
    - 2 \,\, \frac{{r_{ii}^3}}{{h_i^2}}
    + \sum\limits_l {\frac{{{r_{ii}}{r_{ll}}{r_{il}}}}{{{h_i}{h_l}}}} \right).
    \label{eqn:udenom}
\end{align}
The last 4 terms in this expression do not contain $u_{ik}$ explicitly. We will see
how to deal with other terms in this expression whose dominators containing $u_{ik}$
or $u_{il}$. First observe that
\[ \frac{v_{lk}}{u_{ik}} = r_{lk} \left( -1 + \frac{u_i - u_l}{u_i - u_k} \right). \]
Multiplying both sides of this equation by $\theta_{il}$ and using definition of $\theta_{il}$,
we have
\[ \frac{v_{lk} \theta_{il} }{u_{ik}}
    = - r_{lk} \theta_{il} - \frac{r_{lk}}{u_i - u_k} \left(r_{il} + \sum_m r_{im} v_{lm} \right).\]
Plugging  this formula into the first term of expression~\eqref{eqn:udenom}, then using definition
of $\Omega_{ik}$ and equations \eqref{eq:string}, \eqref{eq:theta} and \eqref{eq:g=1 1-pt small}, we obtain the following
formula for the first 4 terms in expression~\eqref{eqn:udenom}:

\begin{align*}
\sum\limits_{\mathop {k,l}\limits_{k \ne i,l \ne i} } {\frac{{{h_k}{v_{lk}}
        {\theta _{il}}}}{{{u_{ik}}h_i^3}}}
- 2 \sum\limits_{k \ne i} {\frac{{{\theta _{ik}}{h_k}}}{{{u_{ik}}h_i^3}}}
+ \frac{1}{{12}}\sum\limits_{\mathop {k,l}\limits_{k \ne i} }
    {\frac{{{r_{ik}}{r_{il}}{h_k}}}{{{u_{ik}}h_i^2{h_l}}}}
- \sum\limits_{k \ne i} {\frac{{{r_{ii}}{r_{ik}}{h_k}}}{{{u_{ik}}h_i^3}}}
=&\sum\limits_{\mathop k\limits_{k \ne i} } {\frac{{{h_k}}}{{h_i^3}}{\Omega _{ik}}}
+ \sum\limits_{\mathop l\limits_{l \ne i} } {\frac{1}{{h_i^2}}{r_{il}}{\theta _{il}}}.
\end{align*}
Similarly, for the next 5 terms in expression~\eqref{eqn:udenom}, we have the following formula
\begin{align*}
 &- 2\sum\limits_{k \ne i} {\frac{{{\theta _{ki}}}}{{{u_{ik}}{h_i}{h_k}}}}
 + \frac{1}{12}\sum\limits_{\mathop {k,l}\limits_{k \ne i} }
    {\frac{{{r_{ik}}{r_{kl}}}}{{{u_{ik}}{h_i}{h_l}}}}
 - 2 \sum\limits_{k \ne i} {\frac{{{r_{kk}}{r_{ik}}}}{{{u_{ik}}{h_i}{h_k}}}}
  + \sum\limits_{\mathop {k,l}\limits_{l \ne i} }
    {\frac{{{r_{ik}}{r_{ll}}{v_{kl}}}}{{{u_{il}}{h_i}{h_l}}}}
  + \sum\limits_{\mathop {k,l}\limits_{k \ne l,l \ne i} }
    {\frac{{{u_{lk}}{r_{ik}}{\theta _{lk}}}}{{{u_{il}}{h_i}{h_l}}}}
 \\=& - \sum\limits_{k \ne i} {\frac{{{\Omega _{ki}}}}{{{h_i}{h_k}}}}
 + \sum\limits_{\mathop l\limits_{l \ne i} } {\frac{{{r_{ll}}}}{{{h_i}{h_l}}}{\theta _{il}}} .
 \end{align*}
Consequently, $\frac{1}{2}P_{ii}^{(2)}+Q_{i}^{(2)}$ can be expressed as a function of
$\{ \Omega_{ij} \}$,  $\{\theta_{jk} \}$, $\{ v_{jk} \}$,  $\{ r_{jk} \}$ and $\{ h_{j} \}$.
The precise form of this function is given by equation $(\ref{Q+P})$ in Appendix~\ref{sec:PQ}.
To show $\frac{1}{2}P_{ii}^{(2)}+Q_{i}^{(2)}=0$, we need to get rid of
$\{ \Omega_{ij} \}$,  $\{\theta_{jk} \}$, $\{ v_{jk} \}$ from this function.

After getting rid of part of $\{\theta_{jk} \}$ terms using equations \eqref{eq:theta}, \eqref{eq:g=1 2-pt small},
Lemma~\ref{lem:Gi} and Remark~\ref{rem:Gi}, equation $(\ref{Q+P})$ has the following form
\begin{align}
&5760\{\frac{1}{2}P_{ii}^{(2)}+Q_{i}^{(2)}\}\nonumber
\\=&  - 5\sum\limits_{\mathop {{\rm{ }}j}\limits_{j \ne i} } {\frac{{{{\widetilde \Omega }_{ij}}}}{{h_i^2}}}   - 100\sum\limits_{k \ne i} {\frac{{{r_{ik}}{\theta _{ik}}}}{{h_i^2}}}
  + 4\sum\limits_{k \ne i} {\frac{{{r_{kk}}{h_i}{\theta _{ik}}}}{{h_k^3}}}  - 6\sum\limits_{k \ne i} {\frac{{h_i^2{r_{ik}}{\theta _{ki}}}}{{h_k^4}}} \nonumber
\\&  + 5\sum\limits_{\mathop {k,l}\limits_{k \ne i} } {\frac{{{r_{kl}}{\theta _{ik}}}}{{{h_i}{h_l}}}} - 2\sum\limits_{\mathop {k,l}\limits_{l \ne i} } {\frac{{{r_{kl}}{h_i}{\theta _{il}}}}{{{h_k}h_l^2}}}+ 40\sum\limits_{\mathop {k,l}\limits_{l \ne i} } {\frac{{{r_{ik}}{v_{lk}}{\theta _{il}}}}{{h_i^2}}}
\nonumber
\\&+\mbox{\{some function of $r$ and $h$\}}
   \label{Q+P jian}
           .\end{align}
              where $\widetilde{\Omega}_{ij}=\Omega_{ij}\{\frac{h_{i}}{h_{j}}-
\frac{h_{j}}{h_{i}}\}$. Obviously, $\widetilde{\Omega}_{ij}=-\widetilde{\Omega}_{ji}$.
In this formula, the most difficult term to deal with is
\[ 40\sum\limits_{\mathop {k,l}\limits_{l \ne i} } {\frac{{{r_{ik}}{v_{lk}}{\theta _{il}}}}{{h_i^2}}}.\]
In order to get rid of this term, we need condition (C3). This is the only place in the
proof of Theorem~\ref{thm:main} where condition (C3) is needed.

\begin{lem} \label{lem:thetarv}
If conditions (C1), (C2), (C3) are satisfied, we have
\begin{align}
&80\sum\limits_{\mathop {{\rm{ }}j}\limits_{j \ne i} } {\sum\limits_k {{\theta _{ij}}{r_{ik}}{v_{jk}}\frac{1}{{h_i^2}}} }
 \nonumber
\\=& - 15\frac{1}{{h_i^2}}\sum\limits_{\mathop j\limits_{j \ne i} } {{{\widetilde \Omega }_{ji}}}  - 5\sum\limits_{\mathop j\limits_{j \ne i} } {{{\widetilde \Omega }_{ji}}\frac{1}{{h_j^2}}}  - 24\sum\limits_{\mathop j\limits_{j \ne i} } {{\theta _{ji}}{r_{jj}}\frac{{{h_i}}}{{h_j^3}}}
 - 400\sum\limits_{\mathop {{\rm{ }}j}\limits_{j \ne i} } {{\theta _{ji}}{r_{ji}}\frac{1}{{h_i^2}}} \nonumber
 \\&+ 22\sum\limits_{\mathop {{\rm{ }}j}\limits_{j \ne i} } {\sum\limits_k {{\theta _{ji}}{r_{jk}}\frac{1}{{{h_k}{h_i}}}} }
 +\{\mbox{some function of $r$ and $h$}\}
 \label{TE der g=2jian}
 \end{align}
for all $i$. The precise form of this formula is given by equation \eqref{lem 2.5} in Appendix B.
\end{lem}
{\bf Proof}:
Condition (C3) is equivalent to
\[ \gwiitwo{\gamma_{\alpha}}|_{\mathcal{H}}=0  \hspace{20pt} \mbox{ and  } \hspace{20pt}
\gwiitwo{\tau_{1}(\gamma_{\alpha})}|_{\mathcal{H}}=0\]
 for all $\alpha$, where $\mathcal{H}$ is the small phase space.
Since
$$  \gwiitwo{ T({\mathcal{E}_i})}
=  \gwiitwo{ \tau_{+}({\mathcal{E}_i})}- \gwii{
\mathcal{E}_{i}\gamma^{\alpha}} \gwiitwo{\gamma_{\alpha}}
$$
and $\mathcal{E}_{i}$ are primary vector fields, we have
\begin{equation}  {\gwiitwo{ T({\mathcal{E}_i})}} |_{\mathcal{H}}=0.
\end{equation}

On the other hand, a formula for $F_2$ was proved in \cite{L07}.
This formula is included in Appendix~\ref{sec:F2} for convenience.
Using this formula and basic properties of the operator $T$ given in equation \eqref{eqn:Tei},
 we can calculate
$ \gwiitwo{T({\mathcal{E}_i})} = T({\mathcal{E}_i}) F_2$.
The formula for $\gwiitwo{T({\mathcal{E}_i})}$ obtained this way is rather complicated. But the
restriction of this formula to the small phase space has a much simpler form.
For example, the restriction of $< \tau_{-}^k(S), \ve_i>$ to the small phase space is $0$
for all $k \geq 1$. Moreover, if conditions (C1) and (C2) are satisfied, we can
use equation \eqref{eq:g=1 2-pt small}, Lemma~\ref{lem:Gi} and Remark~\ref{rem:Gi} to
get rid of many terms containing $\theta_{jk}$. We also observe that the product $\theta_{ij} v_{ij}$
does not contain poles. More precisely
 \begin{align}
 \theta_{ij} v_{ij} = r_{ij} ( r_{ij} + \sum_k r_{ik} v_{jk} )
 \label{theta v}
 \end{align}
by definition of $\theta_{ij}$ and $v_{ij}$. After all these simplifications and
combinations of like terms, the remaining
terms containing $\theta_{ij}$ and $\Omega_{ij}$ are exactly those terms appearing in
equation~\eqref{TE der g=2jian}.
In particular, the left hand side of equation \eqref{TE der g=2jian} comes from
the following two terms which appear in $5760 \gwiitwo{T({\mathcal{E}_i})}|_{\mathcal{H}}$:
 \begin{align}
 &170\sum\limits_{\mathop j\limits_{j \ne i} } {\sum\limits_k {{\theta _{ij}}{r_{jk}}{v_{ik}}\frac{1}{h_{i}^{2}}} }
 + 90\sum\limits_{\mathop j\limits_{j \ne i} } \sum_{k} {\theta _{ji}r_{ik}v_{jk}\frac{1}{h_{i}^{2}}}
  \label{170 90}.
 \end{align}
 For the first term, we can decompose the function $v_{ik}$ into two parts
\[v_{ik}=r_{ik}u_{kj}+r_{ik}u_{ji} \]
and observe that $\theta_{ij} u_{ji}$ no longer has poles.
For the second term in expression \eqref{170 90}, we can replace $\theta_{ji}$ by $-\theta_{ij}$ using
equation~\eqref{eq:theta}. After combining these terms together, expression \eqref{170 90} becomes
\begin{equation}  \label{eqn:32left}
80\sum\limits_{\mathop {j}\limits_{j \ne i} } {\sum\limits_k {{\theta _{ij}}{r_{ik}}{v_{jk}}\frac{1}{{h_i^2}}} }
\end{equation}
 plus a function
which only depends on $\{r_{jk}\}$, $\{h_k\}$, and $\{ v_{jk} \}$.
The expression \eqref{eqn:32left} is exactly the left hand side of equation \eqref{TE der g=2jian}.

We should also note that there are many terms in $\gwiitwo{T({\mathcal{E}_i})}|_{\mathcal{H}}$ which
do not contain $\{ \theta_{jk} \}$, but contain $\{ v_{jk} \}$.
Many of such terms can be get rid of using equation~\eqref{eq:g=1 1-pt small}.
Other terms of this type can be canceled with terms produced in the process of obtaining expression \eqref{eqn:32left}
and terms produced from applying
equation \eqref{theta v} to those terms containing $\{\theta_{jk}v_{jk}\}$.
Eventually,
after straightforward calculations, we
obtain equation \eqref{TE der g=2jian} from
${\gwiitwo{ T({\mathcal{E}_i})}} |_{\mathcal{H}}=0$.
$\Box$

By equations \eqref{Q+P jian}and  \eqref{TE der g=2jian},
we obtain
 \begin{align}
&5760\{P_{ii}^{(2)}+2Q_{i}^{(2)}\}\nonumber
\\=& - 5\frac{1}{{h_i^2}}\sum\limits_{\mathop {{\rm{ }}j}\limits_{j \ne i} } {{{\widetilde \Omega }_{ji}}}  - 5\sum\limits_{\mathop {{\rm{ }}j}\limits_{j \ne i} } {{{\widetilde \Omega }_{ji}}\frac{1}{{h_j^2}}}  + 10\sum\limits_{\mathop {k,l}\limits_{k \ne i} } {\frac{{{r_{kl}}{\theta _{ik}}}}{{{h_i}{h_l}}}}  - 200\sum\limits_{\mathop k\limits_{k \ne i} } {{\theta _{ki}}{r_{ki}}\frac{1}{{h_i^2}}}\nonumber
\\& - 4\sum\limits_{\mathop {k,l}\limits_{l \ne i} } {\frac{{{r_{kl}}{h_i}{\theta _{il}}}}{{{h_k}h_l^2}}}  + 32\sum\limits_{\mathop k\limits_{k \ne i} } {{\theta _{ik}}{r_{kk}}\frac{{{h_i}}}{{h_k^3}}}  + 22\sum\limits_{\mathop j\limits_{j \ne i} } {\sum\limits_k {{\theta _{ji}}{r_{jk}}\frac{1}{{{h_k}{h_i}}}} }  - 12\sum\limits_{k \ne i} {\frac{{h_i^2{r_{ik}}{\theta _{ki}}}}{{h_k^4}}}\nonumber
\\&+\mbox{\{some function of $r$ and $h$\}}\label{Q+Pjian2}
           \end{align}
Here we omit the terms not containing $\{ \theta_{jk} \}$ or $\{ \Omega_{jk} \}$ .
The following lemma will be used to get rid of $\widetilde \Omega_{ij}$ in this formula.

\begin{lem} \label{lem:Omegaij}
Under conditions $(C1)$ and $(C2)$, on the small phase space we have
\begin{align}
0=&{{\widetilde \Omega }_{ij}} + {\theta _{ij}}\{ 20{r_{ij}} + 4{r_{ii}}\frac{{{h_j}}}{{{h_i}}}
- \sum\limits_k {} ({r_{ik}}\frac{{{h_j}}}{{{h_k}}} + {r_{jk}}\frac{{{h_i}}}{{{h_k}}})\}
+ 2\sum\limits_{k \ne i} {{\theta _{ik}}} {r_{jk}}\frac{{{h_j}}}{{{h_i}}}
   \nonumber \\&+\{\mbox{some function of $r$ and $h$}\}
   \label{g=13-pt2jian}
\end{align}
for all $i\neq j$. The complete formula is given by equation (\ref{g=13-pt2}) in Appendix~\ref{sec:PQ}.
\end{lem}
{\bf Proof}: By  equation (\ref{eq:g=1 k-pt}), $\phi_{iij} \mid_{\mathcal{H}} = 0$.
A formula of $\phi_{iij}$ is given in Appendix~\ref{sec:3ptg1}. Using this formula together with
 equation (\ref{eq:g=1 2-pt small}) and Remark~\ref{rem:Gi},
we obtain equation \eqref{g=13-pt2jian}.
$\Box$

Using equation (\ref{g=13-pt2jian}) to get rid of $\Omega_{ij}$ in
 equation (\ref{Q+Pjian2}) and using equation (\ref{eq:g=1 2-pt small}), Lemma~\ref{lem:Gi}
and Remark~\ref{rem:Gi} to simplify the expression, we obtain
\begin{align}
5760\{P_{ii}^{(2)}+2Q_{i}^{(2)}\}
=&32\sum\limits_{\mathop {{\rm{ }}k}\limits_{k \ne i} } {} {\theta _{ik}}{r_{ik}}\frac{1}{{h_i^2}}
    - 16\sum\limits_{\mathop {j,k}\limits_{k \ne i} } {} {\theta _{ik}}{r_{jk}}\frac{{1}}{{{h_i}{h_j}}}
    + 32\sum\limits_{\mathop {{\rm{ }}k}\limits_{k \ne i} }
    {\theta _{ik}}{r_{kk}}\frac{{{h_i}}}{{h_k^3}}  \nonumber \\
& +\mbox{\{some function of $r$ and $h$\}}
\label{eqn:PQ3theta}
\end{align}
The following lemma to needed in dealing with the first three terms on the right hand side of
this equation.
\begin{lem}
 Under conditions $(C1)$ and $(C2)$, on small phase space we have
  \begin{align}
   &4\sum\limits_{\mathop {{\rm{ }}k}\limits_{k \ne i} }{\theta _{ik}}{r_{ik}}\frac{1}{{h_i^2}}
   - 2\sum\limits_{\mathop {j,k}\limits_{k \ne i} } {} {\theta _{ik}}{r_{jk}}\frac{{1}}{{{h_i}{h_j}}} + 4\sum\limits_{\mathop {{\rm{ }}k}\limits_{k \ne i} } {} {\theta _{ik}}{r_{kk}}\frac{{{h_i}}}{{h_k^3}}
   \nonumber
   \\=& \, \{\mbox{some function of $r$ and $h$}\} \label{g=1 3-pt3}
     \end{align}
 for all $i$.  The precise form of this formula is equation (\ref{lem 2.4.2}) in Appendix~\ref{sec:PQ}.
\end{lem}
{\bf Proof:}
We first multiply both sides of equation (\ref{g=1 3-ptijkjian}) by $\frac{1}{h_{k}^{2}}$
and take summation over $j$ and $k$ for distinct $i$, $j$, $k$.
Then by equations (\ref{eq:theta}), (\ref{eq:g=1 2-pt small}), Lemma~\ref{lem:Gi} and Remark~\ref{rem:Gi},
 we obtain equation (\ref{g=1 3-pt3}).$\Box$

Plugging equation \eqref{g=1 3-pt3} into equation \eqref{eqn:PQ3theta},
we can express $\frac{1}{2}P_{ii}^{(2)} + Q_i^{(2)}$ as a function which only depends on
$\{r_{jk} \}$ and $\{h_j\}$.
A straightforward but complicated calculation then shows that this function is identically equal to $0$.
This proves  $\frac{1}{2}P_{ii}^{(2)} + Q_i^{(2)}=0$ for all $i$, and thus finishes the proof of
Theorem~\ref{thm:main}.

\section{Proof of corollary 0.2}
\label{sec:pfcor}

In this section we will prove $G^{(2)}=0$ for the cohomological field theories associated to $ADE$ singularities
by showing that conditions (C1), (C2), (C3) in Theorem~\ref{thm:main} are satisfied.

In fact, conditions (C2) and (C3) just follow from a simple dimension count.
For $ADE$ singularities,  a necessary condition for $\gwig{\tau_{l_{1}}(\gamma_{k_{1}})\cdot\cdot\cdot \tau_{l_{s}}(\gamma_{k_{s}})} \neq 0$ is the following dimension condition
$$2\big((\hat{c}_{W}-3)(1-g)+s\big)=\sum_{i=1}^{s}(2l_{i}+ \deg_{W}\gamma_{k_{i}})$$
where $W$ is a quasi-homogenous polynomial with a singularity of $ADE$ type, $\hat{c}_{W}$ is the central charge of $W$ and
$\deg_{W}\gamma_{i}$ is the $W$-degree of $\gamma_{i}$
(cf: \cite{FJR} for details). For $ADE$ singularities ,
$\hat{c}_{W}<1$, $\deg_{W}\gamma_{i}<2$ for all $i$. Therefore the dimension condition implies that
\[ \gwih{1}{\gamma_{k_{1}} \cdots \gamma_{k_{s}} } =0 \]
and
\begin{equation} \label{eqn:dim}
\gwih{2}{\tau_{m}(\gamma_{k_{1}}) \gamma_{k_{2}} \cdots \gamma_{k_{s}} } =0
\end{equation}
for all $k_i$ and $m \leq 2$. In particular, conditions (C2) and (C3) are satisfied.

We present two different approaches to condition (C1).  The first approach does not involve
canonical coordinates for semisimple Frobenius manifolds. It only uses flat coordinates.
Recall that the genus-2 Mumford relation has the following form (cf. \cite{Ge} and \cite{L07}):
For any vector field $\mathcal{W}$,
\begin{align*}
\gwiitwo{T^{2}(\mathcal{W})}=&\frac{7}{10}\gwiione{\gamma_{\alpha}}
\gwiione{\{\gamma^{\alpha}\circ \mathcal{W}\}}
+\frac{1}{10}\gwiione{\gamma_{\alpha}\{\gua\circ\mathcal{W}\}}-\frac{1}{240}
\gwiione{\mathcal{W}\{\ga\circ\gua\}}
\\&+\frac{13}{240}\gwii{\mathcal{W}\ga\gua\gub}\gwiione{\gb}+\frac{1}{960}
\gwii{\mathcal{W}\ga\gua\gb\gub}.
\end{align*}
By definition of operator $T$, we have
\begin{align*}
T^{2}(\mathcal{W})=\tau_{+}^{2}(\mathcal{W})-\gwii{\mathcal{W}\ga}\tau_{1}(\ga)
-\gwii{T(\mathcal{W})\gua}\ga.
\end{align*}
So
\begin{align*}
\gwiitwo{T^{2}(\mathcal{W})}=\gwiitwo{\tau_{+}^{2}(\mathcal{W})}
-\gwii{\mathcal{W}\ga}\gwiitwo{\tau_{1}(\ga)}
-\gwii{T(\mathcal{W})\gua}\gwiitwo{\ga}.
\end{align*}
By equation~\eqref{eqn:dim}, we have
\begin{align*}
\gwiitwo{T^{2}(\mathcal{W})}\big|_{\mathcal{H}} =0
\end{align*} for any primary vector field $\vw$.
Therefore by Mumford equation and condition (C2), we have
\begin{align} \label{eqn:pfC1}
\gwii{\mathcal{W}\ga\gua\gb\gub}\big|_{\mathcal{H}}=0
\end{align}
for any primary vector field $\mathcal{W}$. This proves
condition (C1) for $ADE$ singularities.

Our second approach to condition (C1) needs \cite[Lemma 2.5]{DLZ}.
During this process, we also found a new proof for \cite[Lemma 2.5]{DLZ}
which is much simpler than the original proof of this lemma in \cite{DLZ}.
The following functions were introduced in \cite{DLZ}:
\[O_{1}=\gwii{\gamma_{\alpha}\gamma_{\alpha'}\gamma_{\beta}
        \gamma_{\beta'}}(M^{-1})^{\alpha\alpha'} (M^{-1})^{\beta\beta'} \]
and
\[ O_{2}=\gwii{\gamma_{\alpha}\gamma_{\beta}\gamma_{\rho}}
    \gwii{\gamma_{1}\gamma_{\alpha'}\gamma_{\beta'}\gamma_{\rho'}}
    (M^{-1})^{\alpha\alpha'}(M^{-1})^{\beta\beta'}(M^{-1})^{\rho\rho'},\]
where entries of the matrix $M$ are defined by
\begin{equation} \label{eqn:M}
 M_{\mu\rho}=\gwii{\gamma_{1}\gamma_{\mu}\gamma_{\rho}}
\end{equation}
 for any $\mu$ and $\rho$.
The following lemma explains the relation between condition (C1) and the condition
\[O_1 - O_2 = \mbox{constant}. \]
\begin{lem} \label{lem:O1O2C1}
For any cohomological field theory,
\begin{align}
\gwiis{\gamma_{\alpha}\gamma^{\alpha}\gamma_{\beta}\gamma^{\beta}} \mid_{s \rightarrow b} \, \, \,\, = \, \,O_{1}-O_{2}
\end{align}
where $\gwiih{g, \mathcal{H}}{\cdots}$ means the restriction of the function $\gwiig{\cdots}$ to the small phase space
$\mathcal{H}$,
and the notation "$f \mid_{s \rightarrow b}$" stands for the function on the big phase space obtained
from a function $f$ on the small phase space via the transformation given by equation~\eqref{eqn:s2b}.
\end{lem}
{\bf Proof:}
The genus-0 constitutive relation has the following form (cf. \cite{DW}):
\[ \gwii{\gamma_{\alpha}\gamma_{\beta}}
=\gwiis{\gamma_{\alpha}\gamma_{\beta}} \mid_{s \rightarrow b}.\]
Taking derivative with respect to $\frac{\partial}{\partial t_{0}^{\rho}}$,
we obtain
$$\gwii{\gamma_{\alpha}\gamma_{\beta}\gamma_{\rho}}
=\gwiis{\gamma_{\alpha}\gamma_{\beta}\gamma^{\mu}} \mid_{s \rightarrow b} M_{\mu\rho},$$
i.e.
\begin{equation} \label{eqn:g0p3const}
\gwiis{\gamma_{\alpha}\gamma_{\beta}\gamma^{\mu}} \mid_{s \rightarrow b}
=\gwii{\gamma_{\alpha}\gamma_{\beta}\gamma_{\rho}}
(M^{-1})^{\rho\mu}.
\end{equation}
Taking derivatives with respect to $\frac{\partial}{\partial t_{0}^{\sigma}}$ again, we obtain
\begin{align*}
\gwii{\gamma_{\alpha}\gamma_{\beta}\gamma_{\rho}\gamma_{\sigma}}
=&\gwiis{\gamma_{\alpha}\gamma_{\beta}\gamma^{\mu}\gamma^{\lambda}}\mid_{s \rightarrow b} M_{\lambda\sigma}M_{\mu\rho}
\\&+
    \gwiis{\gamma_{\alpha}\gamma_{\beta}\gamma^{\mu}}\mid_{s \rightarrow b}
    \gwii{\gamma_{1}\gamma_{\mu}\gamma^{\rho}\gamma^{\sigma}}.
\end{align*}
Multiplying both sides of this equation by $(M^{-1})^{\alpha\rho}(M^{-1})^{\beta\sigma}$, we obtain
\begin{align*}
O_{1}=&\gwiis{\gamma_{\alpha}\gamma_{\beta}\gamma^{\alpha}\gamma^{\beta}}
\mid_{s \rightarrow b}
\\&+\gwiis{\gamma_{\alpha}\gamma_{\beta}\gamma^{\mu}}\mid_{s \rightarrow b}
\gwii{\gamma_{1}\gamma_{\mu}\gamma_{\rho}\gamma_{\sigma}}(M^{-1})^{\alpha\rho}
(M^{-1})^{\beta\sigma}.
\end{align*}
Applying equation~\eqref{eqn:g0p3const} to the second term on the right, we obtain
\[ O_1 = \gwiis{\gamma_{\alpha}\gamma_{\beta}\gamma^{\alpha}\gamma^{\beta}} \mid_{s \rightarrow b} +O_{2}. \]
The lemma is thus proved.
$\Box$

\begin{rem}
For semisimple Frobenius manifolds, it was proved in \cite[Lemma 2.4]{DLZ} that on the small phase space
\[ O_1 - O_2 \,\, = \,\, \sum_{1\leq i<j \leq n}r_{ij}\frac{(h_{i}^{2}+h_{j}^{2})^{2}}{h_{i}^{3}h_{j}^{3}}. \]
After expanding the square in the numerator, it is easy to see that the right hand side of this
equation is precisely the restriction of
$\gwii{\gamma_{\alpha}\gamma^{\alpha} \gamma_{\beta}\gamma^{\beta}}$
to the small phase space due to equations (\ref{eq:string}) and (\ref{eq:O}).
This gives an alternative proof for Lemma~\ref{lem:O1O2C1} for the special case of
semisimple Frobenius manifolds.
\end{rem}

By \cite[Lemma 2.5]{DLZ}, $O_1 - O_2$ is constant
for Frobenius manifolds associated to $ADE$ singularities and $\mathbb{P}^{1}$-orbifold of
$AD$ type. Therefore Lemma~\ref{lem:O1O2C1} implies that condition (C1) also holds for such
Frobenius manifolds.

\begin{rem} \label{rem:pfDLZ2.5}
Note that the proof of \cite[Lemma 2.5]{DLZ} was quite complicated and was done via case by case study.
Lemma~\ref{lem:O1O2C1} and our first proof of equation~\eqref{eqn:pfC1} provide
 a unified and much simpler proof of \cite[Lemma 2.5]{DLZ} for
$ADE$ singularities.
\end{rem}

Since we have verified all three conditions (C1)--(C3), Corollary~\ref{cor:ADEG} just follows from Theorem~\ref{thm:main}.

\section{Further remarks}
\label{sec:rem}

Due to decomposition \eqref{F2}, if $G^{(2)}=0$, the genus-2 potential function $F_2$ can be written
as a linear combination of 16 terms $Q_1$ to $Q_{16}$ as given in \cite[Theorem 1.1]{DLZ}. Each term $Q_i$
can be expressed using flat coordinates and  graphic representations of these terms were also given
in \cite{DLZ}. If in addition, all genus-1 primary invariants are 0, i.e. condition (C2) holds, then
the last 4 terms $Q_{13}$ to $Q_{16}$ are redundant. In this case, $F_2$ can be expressed as a linear combination
of the first 12 terms $Q_1$ to $Q_{12}$. For example, this is true for cohomological field theories
satisfying conditions (C1)--(C3), in particular for $ADE$ singularities.
Recall the following definition from \cite{DLZ}:
\begin{align*}
&{Q_{13}} ={ \gwiione{ {\gamma _\mu }{\gamma _\sigma }}
    \gwii{ {\gamma _1}{\gamma _1}{\gamma _{\mu '}}{\gamma _{\sigma '}}}
    {{({M^{ - 1}})}^{\mu \mu '}}{{({M^{ - 1}})}^{\sigma \sigma '}}},
\\&Q_{14}=\gwiione{\gamma_{1}\gamma_{\alpha}\gamma_{\beta}}
              (M^{-1})^{\alpha\beta},
\\&{Q_{15}} =  { \gwii{ {\gamma _1}{\gamma _\alpha }{\gamma _{\alpha '}}{\gamma _\beta }{ }}
    \gwiione{ {\gamma _1}{\gamma _{\beta '}}}} {({M^{ - 1}})^{\alpha \alpha '}}{({M^{ - 1}})^{\beta \beta '}},
\\&{Q_{16}} =  { \gwiione{ {\gamma _\beta }}  \gwiione{ {\gamma _1}{\gamma _{\beta '}}}} {({M^{ - 1}})^{\beta \beta '}}
             \end{align*}
where matrix $M$ is defined by equation~\eqref{eqn:M}.
The following lemma
describes the precise way to get rid of these terms  under condition (C2).

 \begin{lem} \label{lem:4Q}
  If all genus-1 primary invariants are zero, then
 \begin{align*}
 (i)& \quad  {Q_{16}} - \frac{1}{{24}}{Q_{15}} = 0\\
 (ii)&\quad Q_{ 13}  = \frac{1}{{24}}{Q_3} - \frac{1}{{24}}{Q_4}\\
 (iii)&\quad Q_{ 14} = \frac{1}{{24}}{Q_1} - \frac{1}{{12}}{Q_2} - \frac{1}{{24}}{Q_3} + \frac{1}{{12}}{Q_4}.
 \end{align*}
 \end{lem}
{\bf Proof}: If all the genus-1 primary invariants are zero, then the genus-1 constitutive relation has the form
(cf. \cite{DW}):
 \begin{align*}
     {F_1} = \frac{1}{{24}} {\log \det \{ {\eta ^{ - 1}}M\} }+ \mbox{constant}.
\end{align*}
Note that derivatives of $M^{-1}$ are given by
\[ \frac{\partial }{{\partial t_0^\sigma }}{({M^{ - 1}})^{\alpha \beta }}
     =  -  \gwii{ {\gamma _1}{\gamma _{\alpha '}}{\gamma _{\beta '}}{\gamma _\sigma }}
        ({M^{ - 1}})^{\alpha \alpha '}({M^{ - 1}})^{\beta \beta '}.
\]
Taking derivatives of $F_1$, we get
\begin{align} \label{eqn:dF1}
     \gwiione{ {\gamma _\mu }} =& \frac{1}{{24}} { \gwii{ {\gamma _1}{\gamma _\alpha }{\gamma _\beta }{\gamma _\mu }}{{({M^{ - 1}})}^{\alpha \beta }}},
\end{align}
and
\begin{align} \label{eqn:ddF1}
 \gwiione{ {\gamma _\mu }{\gamma _\sigma }}
  =& - \frac{1}{{24}}{ \gwii{ {\gamma _1}{\gamma _\alpha }{\gamma _\beta }{\gamma _\mu }}
  \gwii{ {\gamma _1}{\gamma _{\alpha '}}{\gamma _{\beta '}}{\gamma _\sigma }}
  {{({M^{ - 1}})}^{\alpha \alpha '}}{{({M^{ - 1}})}^{\beta \beta '}}} \nonumber
  \\& + \frac{1}{{24}}\sum\limits_{\alpha ,\beta }
    { \gwii{ {\gamma _1}{\gamma _\alpha }{\gamma _\beta }{\gamma _\mu }{\gamma _\sigma }}{{({M^{ - 1}})}^{\alpha \beta }}}.
\end{align}
Applying equation~\eqref{eqn:dF1} to the genus-1 1-point function in the definition of $Q_{16}$, we get
$\frac{1}{24} Q_{15}$. This proves part (i). Part (ii) is obtained by applying equation~\eqref{eqn:ddF1} to
the genus-1 2-point function in the definition of $Q_{13}$. Taking derivative of equation~\eqref{eqn:ddF1} with
respect to $\gamma_1$, we get a formula expressing $\gwiione{\gamma_1 \gamma _{\mu} \gamma _{\sigma}}$ in terms of
genus-0 functions, which implies part (iii) after plugging into the definition of $Q_{14}$.
$\Box$

\begin{rem}
The involvement of $Q_{15}$ and $Q_{16}$ in $F_2$ is through the expression
\[
- \frac{7}{{240}}{Q_{15}} + \frac{7}{{10}}{Q_{16}}.
\]
Therefore part (i) of Lemma~\ref{lem:4Q} eliminates both $Q_{15}$ and $Q_{16}$ from $F_{2}$.
\end{rem}

\vspace{40pt}

\appendix
\centerline {\bf \Large Appendix}
\section{The genus-2 G-function}
\label{sec:G2}

In this appendix, we give the precise definition of the genus-2 G-function $G^{(2)}$ following \cite{DLZ}.
Write
\begin{align}
                   G^{(2)}
                   =&\sum_{i}G_{i}^{(2)}(u,u_{x})u_{xx}^{i}
                   +\sum_{i\neq j}G_{ij}^{(2)}(u)\frac{(u_{x}^{j})^{3}}{u_{x}^{i}}
                   \nonumber\\&+\frac{1}{2}\sum_{i,j}P_{ij}^{(2)}(u)u_{x}^{i}u_{x}^{j}
                   +\sum_{i}Q_{i}^{(2)}(u)(u_{x}^{i})^{2}.
                    \label{eqn:g2G}
\end{align}
Let $\gamma_{ij}$ be the rotation coefficient on the small phase space as defined in \cite{D}. Note that
$\gamma_{ii}=0$ which is different from our definition of $r_{ii}$ in Section~\ref{sec:prel}. For $i \neq j$,
$\gamma_{ij}$ is equal to our definition of $r_{ij}$  restricted to the small phase space. Define
\begin{align*}
H_{i}:=\frac{1}{2}\sum_{j\neq i}u_{ij}\gamma_{ij}^{2}
\end{align*}
where $u_{ij}:=u_{i}-u_{j}$.
Then the function $ G_{i}^{(2)}$ can be defined as
\begin{align*}
   G_{i}^{(2)}=G_{i,1}^{(2)}+G_{i,2}^{(2)}
  \end{align*}
with
\begin{align*}
G_{i,1}^{(2)}=&\frac{\partial_{x}h_{i}H_{i}}{60u_{i,x}h_{i}^{3}}
-\frac{7\partial_{i}h_{i}\partial_{x}h_{i}}{5760u_{i,x}h_{i}^{4}}
+\sum_{k}\Bigg(\frac{\gamma_{ik}H_{k}}{120h_{i}h_{k}}\frac{u_{k,x}}{u_{i,x}}
-\frac{\gamma_{ik}\partial_{x}h_{i}}{5760h_{i}^{2}h_{k}u_{i,x}}
-\frac{\gamma_{ik}\partial_{k}h_{k}u_{k,x}}{1152h_{i}h_{k}^{2}u_{i,x}}
\\&+\frac{\partial_{i}\gamma_{ik}h_{k}u_{k,x}}{1920u_{i,x}h_{i}^{3}}
+\frac{\partial_{x}\gamma_{ik}}{5760u_{i,x}h_{i}h_{k}}
+\frac{\partial_{k}\gamma_{ik}u_{k,x}}{2880h_{i}h_{k}u_{i,x}}
-\frac{7\gamma_{ik}^{2}u_{k,x}}{1152h_{i}^{2}u_{i,x}}\Bigg)
-\sum_{k,l}\frac{u_{k,x}h_{k}\gamma_{il}\gamma_{kl}}{1920u_{i,x}h_{i}h_{l}^{2}}
\end{align*}
and
\begin{align*}
G_{i,2}^{(2)}=&-\frac{3\partial_{i}h_{i}H_{i}}{40h_{i}^{3}}+\frac{19(\partial_{i}h_{i})^{2}}{2880h_{i}^{4}}
+\sum_{k}\Bigg(\frac{\gamma_{ik}H_{i}}{120h_{i}h_{k}}+\frac{7\gamma_{ik}H_{k}}{120h_{i}h_{k}}-\frac{4\gamma_{ik}\partial_{i}h_{i}}{5760h_{i}^{2}h_{k}}
-\frac{7\gamma_{ik}\partial_{k}h_{k}}{2880h_{i}h_{k}^{2}}\\&+\frac{\gamma_{ik}\partial_{k}h_{k}}{384h_{i}^{3}}-\frac{\partial_{k}\gamma_{ik}h_{k}}{384h_{i}^{3}}
+\frac{\partial_{i}\gamma_{ik}}{2880h_{i}h_{k}}+\frac{7\partial_{k}\gamma_{ik}}{2880h_{i}h_{k}}+\frac{\gamma_{ik}h_{i}\partial_{k}h_{k}}{2880h_{k}^{4}}
-\frac{19\gamma_{ik}^{2}}{720h_{i}^{2}}+\frac{\gamma_{ik}^{2}}{1440h_{k}^{2}}\Bigg)
\\&-\sum_{k,l}\frac{h_{i}\gamma_{il}\gamma_{kl}}{2880h_{k}h_{l}^{2}}.
\end{align*}
Other functions in equation \eqref{eqn:g2G} are defined in the following way:
   \begin{align*}
G_{ij}^{(2)} =&  - \frac{{\gamma _{ij}^2{H_j}}}{{120h_j^2}} + \frac{{\gamma _{ij}^3}}{{480{h_i}{h_j}}} - \frac{{{\gamma _{ij}}}}{{5760}}(\frac{{{\partial _i}{\gamma _{ij}}}}{{h_i^2}} + \frac{{{\partial _j}{\gamma _{ij}}}}{{h_j^2}}) + \frac{{\gamma _{ij}^2}}{{5760}}(\frac{{{\partial _i}{h_i}}}{{h_i^3}} + \frac{{3{\partial _j}{h_j}}}{{h_j^3}})
\\& + \sum\limits_k {} (\frac{{{\gamma _{ij}}{\gamma _{ik}}{\gamma _{jk}}}}{{5760h_k^2}} + \frac{{\gamma _{ij}^2}}{{5760{h_k}}}(\frac{{{\gamma _{jk}}}}{{{h_j}}} - \frac{{{\gamma _{ik}}}}{{{h_i}}})),
  \end{align*}
   \begin{align*}
P_{ij}^{(2)} =&  - \frac{{2{\gamma _{ij}}{H_i}{H_j}}}{{5{h_i}{h_j}}} + \frac{{{\gamma _{ij}}{\partial _j}{h_j}{H_i}}}{{20{h_i}h_j^2}} + \frac{{{\gamma _{ij}}{h_i}{\partial _j}{h_j}{H_j}}}{{20h_j^4}} - \frac{{19\gamma _{ij}^2{H_j}}}{{30h_j^2}} - \frac{{{\partial _i}{\gamma _{ij}}{H_j}}}{{60{h_i}{h_j}}}
\\& + \frac{{41\gamma _{ij}^3}}{{240{h_i}{h_j}}} - \frac{{41{\gamma _{ij}}{\partial _i}{\gamma _{ij}}}}{{1440h_i^2}} + \frac{{{\partial _i}{\gamma _{ij}}{\partial _j}{h_j}}}{{1440{h_i}h_j^2}} + \frac{{79\gamma _{ij}^2{\partial _j}{h_j}}}{{1440h_j^3}} - \frac{{{\gamma _{ij}}{\partial _i}{h_i}{\partial _j}{h_j}}}{{720h_i^2h_j^2}} - \frac{{{\gamma _{ij}}{h_i}{{({\partial _j}{h_j})}^2}}}{{288h_j^5}}
\\& + \sum\limits_k {} (\frac{{{\gamma _{ij}}{\gamma _{ik}}{H_j}}}{{60{h_j}{h_k}}} - \frac{{{\gamma _{ik}}{\gamma _{jk}}{h_i}{h_j}{H_k}}}{{30h_k^4}} - \frac{{{\gamma _{ij}}{\gamma _{jk}}{h_i}{H_j}}}{{60h_j^2{h_k}}} + \frac{{{\gamma _{ik}}{\gamma _{jk}}{h_i}{H_j}}}{{60{h_j}h_k^2}} - \frac{{7{\gamma _{ij}}{\gamma _{jk}}{h_i}{H_k}}}{{60h_j^2{h_k}}}
\\& - \frac{{{\gamma _{ij}}{\gamma _{ik}}{\partial _j}{h_j}}}{{720h_j^2{h_k}}} + \frac{{{\gamma _{ij}}{\gamma _{jk}}{h_i}{\partial _j}{h_j}}}{{240h_j^3{h_k}}} - \frac{{{\gamma _{ik}}{\gamma _{jk}}{h_i}{\partial _j}{h_j}}}{{1440h_j^2h_k^2}} + \frac{{{\gamma _{ij}}{\gamma _{jk}}{h_i}{\partial _k}{h_k}}}{{720h_k^4}} + \frac{{{\gamma _{ik}}{\gamma _{jk}}{h_i}{h_j}{\partial _k}{h_k}}}{{288h_k^5}}
\\& + \frac{{{\gamma _{jk}}{\partial _i}{\gamma _{ij}}}}{{1440{h_i}{h_k}}} - \frac{{{h_j}{h_k}{\gamma _{ij}}{\partial _i}{\gamma _{ik}}}}{{360h_i^4}} - \frac{{{h_j}(3{\gamma _{ik}}{\partial _i}{\gamma _{ij}} + 2{\gamma _{ij}}{\partial _i}{\gamma _{ik}})}}{{1440h_i^2{h_k}}} - \frac{{7{h_j}{\gamma _{ij}}{\partial _k}(h_k^{ - 1}{\gamma _{ik}})}}{{1440h_i^2}}
\\& - \frac{{{h_i}{h_j}{\gamma _{ik}}{\partial _k}{\gamma _{jk}}}}{{480h_k^4}} + \frac{{\gamma _{ij}^2{\gamma _{jk}}}}{{120{h_j}{h_k}}} + \frac{{7{h_i}{\gamma _{ij}}\gamma _{jk}^2}}{{160h_j^3}} + \frac{{11{\gamma _{ij}}{\gamma _{ik}}{\gamma _{jk}}}}{{2880h_k^2}} + \frac{{{h_j}\gamma _{ik}^2{\gamma _{jk}}}}{{96h_k^3}})
\\& + \sum\limits_{k,l} {} (\frac{{{h_i}{h_j}{\gamma _{il}}{\gamma _{jl}}}}{{720{h_k}h_l^2}}(\frac{{{\gamma _{kl}}}}{{{h_l}}} - \frac{{{\gamma _{jk}}}}{{2{h_j}}}) - \frac{{{h_i}{\gamma _{ij}}{\gamma _{jl}}{\gamma _{kl}}}}{{720{h_k}h_l^2}}),
  \end{align*}
and
 \begin{align*}
Q_i^{(2)} =& \frac{{4H_i^3}}{{5h_i^2}} - \frac{{7{\partial _i}{h_i}H_i^2}}{{10h_i^3}} + \frac{{7{{({\partial _i}{h_i})}^2}{H_i}}}{{48h_i^4}} - \frac{{{{({\partial _i}{h_i})}^3}}}{{120h_i^5}} + \sum\limits_k {} (\frac{{7{\gamma _{ik}}{H_i}{H_k}}}{{10{h_i}{h_k}}} - \frac{{{\gamma _{ik}}{\partial _i}{h_i}{H_i}}}{{120h_i^2{h_k}}}
\\& + \frac{{7{\partial _k}(h_k^{ - 1}{\gamma _{ik}}){H_i}}}{{240{h_i}}} - \frac{{7{\gamma _{ik}}{\partial _i}{h_i}{H_k}}}{{80h_i^2{h_k}}} + \frac{{{\gamma _{ik}}{H_k}}}{{576{u_{ik}}{h_i}{h_k}}} + \frac{{(2{H_i} + 7{H_k}){\partial _i}{\gamma _{ik}}}}{{240{h_i}{h_k}}}
\\& + \frac{{{\gamma _{ik}}{h_k}{H_i}}}{{576{u_{ik}}h_i^3}} - \frac{{31\gamma _{ik}^2{H_i}}}{{144h_i^2}} + \frac{{{\gamma _{ik}}{{({\partial _i}{h_i})}^2}}}{{720h_i^3{h_k}}} + \frac{{253\gamma _{ik}^2{\partial _i}{h_i}}}{{5760h_i^3}} - \frac{{{\partial _i}{\gamma _{ik}}{\partial _i}{h_i}}}{{960h_i^2{h_k}}} - \frac{{\gamma _{ik}^2{\partial _k}{h_k}}}{{2880h_k^3}}
\\& - \frac{{7{\partial _k}(h_k^{ - 1}{\gamma _{ik}}){\partial _i}{h_i}}}{{1920h_i^2}} - \frac{{7{\partial _i}{\gamma _{ik}}{\partial _k}{h_k}}}{{5760{h_i}h_k^2}} - \frac{{41{\partial _i}{\gamma _{ik}}{\partial _i}{h_i}{h_k}}}{{5760h_i^4}} + \frac{{{\partial _i}({h_i}{\gamma _{ik}}){\partial _k}{h_k}}}{{2880h_k^4}}
\\& - \frac{{113{\gamma _{ik}}{\partial _i}{\gamma _{ik}}}}{{5760h_i^2}} + \frac{{(3{\partial _i}{\gamma _{ik}} + {\partial _k}{\gamma _{ik}}){\gamma _{ik}}}}{{1440h_k^2}} - \frac{{{\partial _i}{\gamma _{ik}}{h_k}}}{{576{u_{ik}}h_i^3}} - \frac{{{\partial _k}{\gamma _{ik}}}}{{576{u_{ik}}{h_i}{h_k}}} - \frac{{\gamma _{ik}^3}}{{240{h_i}{h_k}}})
\\& + \sum\limits_{k,l} {} ( - \frac{{{\gamma _{kl}}{\partial _i}({h_i}{\gamma _{il}})}}{{2880{h_k}h_l^2}} + \frac{{\gamma _{il}^2{\gamma _{kl}}}}{{2880{h_k}{h_l}}} - \frac{{{\gamma _{ik}}\gamma _{il}^2}}{{240{h_i}{h_k}}} - \frac{{{\gamma _{kl}}{\partial _i}{\gamma _{ik}}}}{{2880{h_i}{h_l}}} + \frac{{{u_{lk}}{\gamma _{ik}}{\partial _l}{\gamma _{kl}}}}{{1152{u_{il}}{h_i}{h_l}}}
\\& + \frac{{{u_{kl}}{\gamma _{ik}}{\gamma _{kl}}{\partial _i}{\gamma _{il}}}}{{144h_i^2}} + \frac{{{h_l}{\gamma _{ik}}{\partial _i}{\gamma _{il}}}}{{144h_i^2{h_k}}} + \frac{{{h_k}{u_{kl}}{\gamma _{kl}}{\partial _i}{\gamma _{il}}}}{{1152{u_{ik}}h_i^3}} + \frac{{{h_l}{u_{ik}}\gamma _{ik}^2{\partial _i}{\gamma _{il}}}}{{40h_i^3}}).
                          \end{align*}

In these expressions, all summations are taken over the ranges of indices where the denominators do not vanish.
The symbol $\partial_x$ means taking derivative with respect to $x=t_0^1$ after transformation
given by equation~\eqref{eqn:s2b} and $\partial_i$ means taking derivative
with respect to $u_i$.

\section{Explicit formulas of some equations appeared in the proof of the main theorem}
\label{sec:PQ}

For clarity of presentations, we omitted terms only depending on $\{ r_{ij} \}$ and $\{h_i\}$ in many formulas
in sections \ref{sec:P=0} and \ref{sec:P+Q=0}. We give precise forms of such formulas in this section for
readers' convenience.

First, after getting rid of functions $\{H_k\}$ and partial derivatives of $\{r_{kl}\}$ and $\{h_k\}$
in the definition of $P_{ij}^{(2)}$, we obtain
the following formula for $P_{ij}^{(2)}$:
\begin{align}
P_{ij}^{(2)}=
&-\frac{1}{480}\sum_{k\neq j}\frac{h_{i}h_{j}r_{ik}\theta_{kj}}{h_{k}^{4}}-\frac{41r_{ij}\theta_{ij}}{1440h_{i}^{2}}
+\frac{r_{ii}h_{j}\theta_{ij}}{240h_{i}^{3}}
-\frac{1}{480}\sum_{k}\frac{r_{ik}h_{j}\theta_{ij}}{h_{i}^{2}h_{k}}
\nonumber
\\&-\frac{1}{360}\sum_{k\neq
i}\frac{r_{ij}h_{j}h_{k}\theta_{ik}}{h_{i}^{4}}
-\frac{1}{288}\sum_{k\neq
i}\frac{r_{ij}h_{j}\theta_{ki}}{h_{i}^{2}h_{k}}+\frac{41r_{ij}^{3}}{240h_{i}h_{j}}
-\frac{13r_{ij}r_{jj}^{2}h_{i}}{240h_{j}^{3}}
-\frac{r_{ij}r_{ii}^{2}h_{j}}{240h_{i}^{3}}
\nonumber\\&-\frac{93r_{ii}r_{ij}^{2}}{2880h_{i}^{2}}
+\frac{93r_{jj}r_{ij}^{2}}{2880h_{j}^{2}}
-\frac{1}{720}\sum_{k}\frac{r_{ij}r_{ik}^{2}h_{j}}{h_{i}^{3}}
+\frac{1}{1440}\sum_{k}\frac{r_{ij}r_{jj}r_{ik}}{h_{j}h_{k}}
-\frac{13}{720}\sum_{k}\frac{r_{ij}^{2}r_{jk}}{h_{j}h_{k}}
\nonumber\\&-\frac{1}{1440}\sum_{k}\frac{r_{ii}r_{ij}r_{jk}}{h_{i}h_{k}}
+\frac{19}{1440}\sum_{k}\frac{r_{ij}r_{jj}r_{jk}h_{i}}{h_{j}^{2}h_{k}}
-\frac{1}{288}\sum_{k}\frac{r_{ij}r_{ii}r_{ik}h_{j}}{h_{i}^{2}h_{k}}
+\frac{7}{160}\sum_{k}\frac{h_{i}r_{ij}r_{jk}^{2}}{h_{j}^{3}}
\nonumber\\&+\frac{11}{2880}\sum_{k}\frac{r_{ij}r_{ik}r_{jk}}{h_{k}^{2}}
+\frac{1}{96}\sum_{k}\frac{r_{ik}^{2}r_{jk}h_{j}}{h_{k}^{3}}
+\frac{1}{360}\sum_{k}\frac{r_{ij}r_{jk}r_{kk}h_{i}}{h_{k}^{3}}
+\frac{1}{720}\sum_{k,l}\frac{r_{ij}r_{il}r_{kl}h_{j}}{h_{i}^{2}h_{k}}
\nonumber\\&-\frac{1}{1440}\sum_{k,l}\frac{r_{ij}r_{jk}r_{jl}h_{i}}{h_{j}h_{k}h_{l}}
-\frac{7}{1440}\sum_{k,l}\frac{r_{ij}r_{jk}r_{kl}h_{i}}{h_{j}^{2}h_{l}}
-\frac{1}{720}\sum_{k,l}\frac{r_{ij}r_{jl}r_{kl}h_{i}}{h_{k}h_{l}^{2}}.
\label{eqn:Pijc}
\end{align}
Similarly, after getting rid of functions $\{H_k\}$, partial derivatives of $\{r_{kl}\}$ and $\{h_k\}$,
and functions $\{ u_k - u_l \}$ explicitly involved
in the definition of $\frac{1}{2}P_{ii}^{(2)}+Q_{i}^{(2)}$, we obtain
the following formula:
\begin{align}
&5760\{\frac{1}{2}P_{ii}^{(2)}+Q_{i}^{(2)}\}\nonumber
\\=&\{  - 5\sum\limits_{\mathop {{\rm{ }}j}\limits_{j \ne i} } {\frac{{{{\widetilde \Omega }_{ij}}}}{{h_i^2}}}  + 5\sum\limits_{\mathop {k,l}\limits_{k \ne i} } {\frac{{{r_{kl}}{\theta _{ik}}}}{{{h_i}{h_l}}}}  - 108\sum\limits_{k \ne i} {\frac{{{r_{ik}}{\theta _{ik}}}}{{h_i^2}}}  -8\sum\limits_{k \ne i} {\frac{{{r_{ik}}{\theta _{ki}}}}{{h_k^2}}}
- 2\sum\limits_{\mathop {k,l}\limits_{l \ne i} } {\frac{{{r_{kl}}{h_i}{\theta _{il}}}}{{{h_k}h_l^2}}} \nonumber
\\&  + 40\sum\limits_{\mathop {k,l}\limits_{l \ne i} } {\frac{{{r_{ik}}{v_{lk}}{\theta _{il}}}}{{h_i^2}}}  + 4\sum\limits_{k \ne i} {\frac{{{r_{kk}}{h_i}{\theta _{ik}}}}{{h_k^3}}}  - 6\sum\limits_{k \ne i} {\frac{{h_i^2{r_{ik}}{\theta _{ki}}}}{{h_k^4}}} \}
  \nonumber
 \\
& + \{  - 40\sum\limits_{k \ne i} {\frac{{{r_{ii}}{h_k}{\theta _{ik}}}}{{h_i^3}}}  + 5\sum\limits_{\mathop {k,l}\limits_{k \ne i} } {\frac{{{r_{il}}{\theta _{ki}}}}{{{h_l}{h_k}}}}  + 16\sum\limits_{\mathop {k,l}\limits_{l \ne i} } {\frac{{{r_{ik}}{h_l}{\theta _{il}}}}{{h_i^2{h_k}}}}  - 10\sum\limits_{k \ne i} {\frac{{{r_{ii}}{\theta _{ki}}}}{{{h_i}{h_k}}}}\}
            \nonumber
           \\& - 240\frac{{r_{ii}^3}}{{h_i^2}} + 210\sum\limits_k {\frac{{r_{ik}^2{r_{ii}}}}{{h_i^2}}}  + 8\sum\limits_k {\frac{{{r_{kk}}{r_{ii}}{r_{ik}}{h_i}}}{{h_k^3}}}  - 24\sum\limits_k {\frac{{r_{ik}^3}}{{{h_i}{h_k}}}}  + 12\sum\limits_k {\frac{{{r_{ii}}r_{ik}^2}}{{h_k^2}}}
            \nonumber
          \\& + 118\sum\limits_k {\frac{{r_{ii}^2{r_{ik}}}}{{{h_i}{h_k}}}}  + 30\sum\limits_k {\frac{{r_{ik}^3{h_i}}}{{h_k^3}}}  - 14\sum\limits_{k,l} {\frac{{{r_{ii}}{r_{ik}}{r_{il}}}}{{{h_k}{h_l}}}}  - 10\sum\limits_{k,l} {\frac{{{r_{ik}}{r_{ii}}{r_{kl}}}}{{{h_i}{h_l}}}}  - 4\sum\limits_{k,l} {\frac{{{r_{kl}}{r_{ii}}{r_{il}}{h_i}}}{{{h_k}h_l^2}}}
           \nonumber
           \\& -12\sum_{k,l}\frac{r_{ik}r_{il}r_{kl}}{h_{k}^{2}}+ 2\sum\limits_{k,l} {\frac{{r_{il}^2{r_{kl}}}}{{{h_k}{h_l}}}}  - \frac{{221}}{3}\sum\limits_{k,l} {\frac{{r_{ik}^2{r_{il}}}}{{{h_i}{h_l}}}}  + \frac{1}{3}\sum\limits_{j,k,l} {\frac{{{r_{ik}}{r_{il}}{r_{ij}}{h_i}}}{{{h_j}{h_k}{h_l}}}}  + 5\sum\limits_{j,k,l} {\frac{{{r_{ik}}{r_{ij}}{r_{kl}}}}{{{h_j}{h_l}}}}  \label{Q+P}
           .\end{align}
 The precise formula for equation (\ref{g=1 3-pt jian}) is the following
\begin{align}
&\theta _{ik}r_{jk}\frac{h_j}{h_i} + \theta _{ki}r_{ij}\frac{h_j}{h_k}+ \theta _{jk}r_{ik}\frac{h_i}{h_j}
 + \theta _{kj}r_{ij}\frac{h_i}{h_k}+ \theta _{ij}r_{jk}\frac{h_k}{h_i}
+ {\theta _{ji}}{r_{ik}}\frac{{{h_k}}}{{{h_j}}}
  \nonumber \\
=&
6 r_{ij}^{2}r_{ik}\frac{h_k}{h_i} + 6 r_{ij}^{2}r_{jk}\frac{h_k}{h_j}
+ 6 r_{ik}^2{r_{ij}}\frac{{{h_j}}}{{{h_i}}}
 + 6 r_{jk}^{2}r_{ij}\frac{h_i}{h_j} + 6 r_{ik}^{2}r_{jk}\frac{h_j}{h_k}
 + 6 r_{jk}^{2}r_{ik}\frac{h_i}{h_k}
                             \nonumber\\&
 + 12 r_{ij}r_{ik}r_{jk}
                            \nonumber
\\&+ \sum\limits_l \{  r_{il}r_{jk}r_{kl}\frac{h_{i}h_{j}}{h_{l}^{2}}
            + r_{ik}r_{jl}r_{kl}\frac{h_{i}h_{j}}{h_{l}^{2}}
            + r_{ik}r_{il}r_{jl}\frac{h_{j}h_{k}}{h_{l}^{2}}
            + r_{ij}r_{il}r_{kl}\frac{{h_{j}h_{k}}}{h_{l}^{2}} \nonumber
\\&\quad\quad\quad + r_{ij}r_{jl}r_{kl}\frac{h_{i}h_{k}}{h_{l}^{2}}
            + r_{il}r_{jk}r_{jl}\frac{h_{i}h_{k}}{h_{l}^{2}}
            - r_{il}r_{jl}r_{kl}\frac{h_{i}h_{j}h_{k}}{h_{l}^{3}}
            - r_{ik}r_{jk}r_{kl}\frac{h_{i}h_{j}}{h_{k}h_{l}} \nonumber
\\& \quad\quad\quad- r_{ij}r_{ik}r_{il}\frac{h_{j}h_{k}}{h_{i}h_{l}}
            - r_{ij}r_{jk}r_{jl}\frac{h_{i}h_{k}}{h_{j}h_{l}}
            - r_{ij}r_{jk}r_{kl}\frac{h_{i}}{h_{l}}
            - r_{ik}r_{jk}r_{jl}\frac{h_{i}}{h_{l}} \nonumber
\\&\quad\quad\quad- r_{ik}r_{jk}r_{il}\frac{h_{j}}{h_{l}}
            - r_{ij}r_{ik}r_{kl}\frac{h_{j}}{h_{l}}
            - r_{ij}r_{il}r_{jk}\frac{h_{k}}{h_{l}}
            - r_{ij}r_{ik}r_{jl}\frac{h_{k}}{h_{l}}\}
                                \label{g=1 3-ptijkjian}
\end{align}
for any $i\neq j\neq k$.

The precise formula for equation (\ref{TE der g=2jian}) is the following:
                     \begin{align}
    &80\sum\limits_{\mathop j\limits_{j \ne i} } {\sum\limits_k {{\theta _{ij}}{r_{ik}}{v_{jk}}\frac{1}{{h_i^2}}} }
    \nonumber\\=&-15\frac{1}{h_i^2}\sum\limits_{\mathop j\limits_{j \ne i} }
    \widetilde{\Omega}_{ji}  - 5\sum\limits_{\mathop j\limits_{j \ne i} }
    \widetilde{\Omega}_{ji} \frac{1}{{h_j^2}}
      - 24\sum\limits_{\mathop {{\rm{ }}j}\limits_{j \ne i} } {{\theta _{ji}}{r_{jj}}\frac{{{h_i}}}{{h_j^3}}}
       - 400\sum\limits_{\mathop j\limits_{j \ne i} } {{\theta _{ji}}{r_{ji}}\frac{1}{{h_i^2}}}  + 22\sum\limits_{\mathop j\limits_{j \ne i} } {\sum\limits_k {{\theta _{ji}}{r_{jk}}\frac{1}{{{h_k}{h_i}}}} }
       \nonumber\\& + 288r_{ii}^{3}\frac{1}{h_{i}^{2}}
          + 792\sum_{j} r_{ii}r_{ij}^{2}\frac{1}{h_{i}^{2}}
           + 180\sum_{j} r_{ii}r_{ij}^{2}\frac{1}{h_{j}^{2}}
          + 144\sum_{j} r_{ij}r_{ii}^{2}\frac{1}{h_{i}h_{j}}
          - 72 \sum_{j} {r_{jj}}{r_{ij}}{r_{ii}}\frac{{{h_i}}}{{h_j^3}}
         \nonumber \\&- 960\sum_{j} r_{ij}^{3}\frac{1}{{h_i}{h_j}}
         - \frac{{284}}{3}  \sum_{j,k}  r_{ik} r_{ij}^{2}\frac{1}{h_{i}h_{k}}
            + 144\sum_{j,k} r_{jk}r_{ij}^{2}\frac{1}{h_{j}h_{k}}
             - 60 \sum_{j,k} r_{ij}r_{ii}r_{ik}\frac{1}{h_{j}h_{k}}
                 \nonumber \\&- 67\sum_{j,k}  r_{ij}^{2}r_{ik}\frac{h_{i}}{h_{k}h_{j}^{2}}
               - 24\sum_{j,k}  r_{kk}r_{jk}r_{ij}\frac{h_{i}}{h_{k}^{3}}
                +36 \sum_{j,k} r_{jk}r_{ij}r_{ii}\frac{h_{i}}{h_{k}h_{j}^{2}}
                + 26\sum_{j,k} r_{jj}r_{ij}r_{ik}\frac{h_{i}^{2}}{h_{k}h_{j}^{3}}
                \nonumber \\&- 400\sum_{j,k} r_{ij}r_{ik}r_{jk}\frac{1}{h_{i}^{2}}
                + 12 \sum_{j,k,l} r_{lj}r_{ij}r_{lk}\frac{h_{i}}{h_{k}h_{j}^{2}}
                + 12 \sum_{j,k,l} r_{li}r_{kj}r_{jl}\frac{h_{i}}{h_{k}h_{j}^{2}}
                - 12\sum_{k,j,l} r_{li}r_{lk}r_{lj}\frac{h_{i}}{h_{l}h_{j}h_{k}}
                 \nonumber  \\&+ 22\sum_{j,k,l} r_{jk}r_{ik}r_{jl}\frac{1}{h_{i}h_{l}}
                 - 12\sum_{j,k,l} r_{ij}r_{il}r_{jk}\frac{1}{h_{k}h_{l}}
                 -13\sum_{j,k,l} r_{ij}r_{ik}r_{jl}\frac{h_{i}^{2}}{h_{j}^{2}h_{k}h_{l}}
                 \nonumber
               \\& + \frac{37}{3}  \sum_{j,k,l} r_{ij}r_{ik}r_{il}\frac{h_{i}}{h_{j}h_{k}h_{l}}
                    \label{lem 2.5}
                     \end{align}
for any $i$.

The precise formula for equation (\ref{g=13-pt2jian}) is the following
       \begin{align}
0=&{{\widetilde \Omega }_{ij}} + {\theta _{ij}}\{ 20{r_{ij}} + 4{r_{ii}}\frac{{{h_j}}}{{{h_i}}} - \sum_{k} ({r_{ik}}\frac{{{h_j}}}{{{h_k}}} + {r_{jk}}\frac{{{h_i}}}{{{h_k}}})\}  + 2\sum_{k \ne i} {{\theta _{ik}}} {r_{jk}}\frac{{{h_j}}}{{{h_i}}}
   \nonumber \\&+ 24r_{ii}^2{r_{ij}}\frac{{{h_j}}}{{{h_i}}} + 24r_{ij}^3\frac{{{h_i}}}{{{h_j}}}
    - 12\sum_{k}  r_{ik}^2{r_{jk}}\frac{{{h_j}}}{{{h_k}}}
   + 4 \sum_{k} {r_{ij}}{{r_{ik}}{r_{jk}}\frac{{h_i^2}}{{h_k^2}}}   - 4 \sum_{k} {{r_{ij}^{2}}{r_{ik}}\frac{{{h_i}}}{{{h_k}}}}
     \nonumber\\& - 4 \sum_{k}  {r_{ii}}{r_{ij}}{r_{ik}}\frac{{{h_j}}}{{{h_k}}}  - 3 \sum_{k} {{r_{ij}}r_{ik}^2\frac{{{h_i}{h_j}}}{{h_k^2}}}
     - 4\sum_{k} {r_{ik}^2{r_{ij}}\frac{{{h_j}}}{{{h_i}}}}  + 2\sum_{k} {{r_{ij}}{r_{kk}}{r_{ik}}\frac{{h_i^2{h_j}}}{{h_k^3}}}
      \nonumber\\&- 2 \sum_{k}   r_{ij}^2{r_{jk}}\frac{{h_i^2}}{{{h_k}{h_j}}}
    - 2 \sum_{k}r_{ik}^2{r_{jk}}\frac{{h_i^2{h_j}}}{{h_k^3}}
     -  \sum_{k,l} {\frac{{h_i^2{h_j}}}{{{h_k}h_l^2}}{r_{ij}}{r_{kl}}{r_{li}}}   +  \sum_{k,l} {\frac{{{h_i}{h_j}}}{{{h_k}{h_l}}}{r_{ij}}{r_{ik}}{r_{li}}}
      \nonumber\\&- \sum_{k,l} {{r_{kl}}{r_{jl}}{r_{ik}}\frac{{{h_i}{h_j}}}{{h_l^2}}}
   - \sum_{k,l} {{r_{il}}{r_{kl}}{r_{jk}}\frac{{{h_i}{h_j}}}{{h_l^2}}}
  + \sum_{k,l} {{r_{ik}}{r_{il}}{r_{jk}}\frac{{{h_j}}}{{{h_l}}}}   + \sum_{k,l} {{r_{ik}}{r_{kl}}{r_{jk}}\frac{{{h_i}{h_j}}}{{{h_k}{h_l}}}}
   \label{g=13-pt2}
    \end{align}
 for any $i\neq j$.

The precise formula for equation (\ref{g=1 3-pt3}) is the following
\begin{align}
   &4\sum\limits_{\mathop {{\rm{ }}k}\limits_{k \ne i} } {\theta _{ik}}{r_{ik}}\frac{1}{{h_i^2}}
   - 2\sum\limits_{\mathop {j,k}\limits_{k \ne i} } {} {\theta _{ik}}{r_{jk}}\frac{{1}}{{{h_i}{h_j}}} + 4\sum\limits_{\mathop {{\rm{ }}k}\limits_{k \ne i} } {} {\theta _{ik}}{r_{kk}}\frac{{{h_i}}}{{h_k^3}}
   \nonumber
   \\=&- 96 r_{ii}^{3}\frac{1}{{h_{i}^{2}}}
 +48\sum\limits_j r_{ij}^{3} \frac{1}{{{h_i}{h_j}}}   + 44\sum\limits_k {} r_{ii}^{2} {r_{ik}}\frac{1}{{{h_i}{h_k}}}
      + 12\sum\limits_k { r_{ik}^2{r_{ii}}\frac{1}{{h_k^2}} }  - 8\sum\limits_k {r_{ii}}{r_{kk}}{r_{ik}}\frac{{{h_i}}}{{h_k^3}}
       \nonumber\\&+ 20 \sum\limits_k r_{ik}^{2}{r_{ii}}\frac{1}{{h_i^2}}   + 28\sum\limits_k { r_{ik}^3\frac{{{h_i}}}{{h_k^3}} }
- 12\sum\limits_{j,k} { {r_{ij}}{r_{ik}}{r_{jk}}\frac{1}{{h_k^2}} }  - 10  \sum\limits_{j,k} { r_{ij}^2{r_{ik}}\frac{1}{{{h_i}{h_k}}} }
  \nonumber\\& - 18\sum\limits_{j,k} {} r_{ij}^2{r_{jk}}\frac{1}{{{h_j}{h_k}}}
 - 8\sum\limits_{j,k}   {r_{ii}}{r_{ij}}{r_{ik}}\frac{1}{{{h_j}{h_k}}} - 7\sum\limits_{j,k} { {r_{ij}}r_{ik}^2}{\frac{{{h_i}}}{{h_k^2{h_j}}}  }  + 4\sum\limits_{j,k} {r_{ii}}{r_{jk}}{r_{ij}}{\frac{{{h_i}}}{{{h_k}h_j^2}}}
 \nonumber\\
 & + 2\sum\limits_{j,k} { {{r_{ik}}{r_{kk}}{r_{ij}}} \frac{{h_i^2}}{{h_k^3{h_j}}} }  + 4\sum\limits_{j,k}   {r_{ik}}{r_{ij}}{r_{kj}}\frac{{h_i^2}}{{h_k^2h_j^2}}  - 4\sum\limits_{j,k}   {r_{ik}^{2}}{r_{jk}}\frac{{h_i^2}}{{h_k^3{h_j}}}
 +2\sum\limits_{j,k,l}   {r_{ij}}{r_{ik}}{r_{jl}}\frac{1}{{{h_k}{h_l}}}
 \nonumber
\\&- 2\sum\limits_{j,k,l} {r_{ik}}{r_{jk}}{r_{jl}}\frac{{{h_i}}}{{{h_l}h_k^2}}  - 2\sum\limits_{j,k,l} {} {{r_{jk}}{r_{kl}}{r_{ij}}\frac{{{h_i}}}{{h_k^2{h_l}}}}
 + 2\sum\limits_{j,k,l}   {r_{ij}}{r_{jk}}{r_{jl}}\frac{{{h_i}}}{{{h_j}{h_k}{h_l}}}  + \sum\limits_{j,k,l} { {{r_{ij}}{r_{il}}{r_{ik}}}\frac{{{h_i}}}{{{h_k}{h_j}{h_l}}}  }
    \nonumber \\&- \sum\limits_{j,k,l} { {{r_{il}}{r_{jk}}{r_{ij}}} \frac{{h_i^2}}{{{h_k}h_j^2{h_l}}} } \label{lem 2.4.2}
     \end{align}
for any $i$.

\section{The genus-2 generating function}
\label{sec:F2}

In \cite{L07}, the first author of this paper proved a formula for
the generating function $F_2$ for genus-2 Gromov-Witten invariants of compact symplectic manifolds with semisimple
quantum cohomology. Since the proof of this formula only used universal equations which were derived from splitting
principle and relations in tautological
rings of  moduli spaces of curves, this formula also holds for semisimple cohomological field theory.
Since this formula is needed in the proof of Lemma~\ref{lem:thetarv}, we include this formula here for readers' convenience.

 \begin{align*}
  &5760F_{2}
  \\=&- 5 \sum\limits_i {\frac{1}{{g_i^2}} < \tau _ - ^3(S),\mathcal{E}_{i} > }
+  \sum\limits_i {\sum\limits_{\mathop j\limits_{j \ne i} } {5{\Omega _{ij}}(\frac{1}{{{g_j}}}\sqrt {\frac{{{g_i}}}{{{g_j}}}}  - \frac{1}{{\sqrt {{g_i}{g_j}} }})} }
\\& +\sum\limits_i {\frac{1}{{{g_i}}} < \tau _ - ^2(S),\mathcal{E}_{i} > \{ 24{r_{ii}}\frac{1}{{{g_i}}}+  \sum\limits_j (5{r_{ij}}\frac{1}{{{g_j}}}\sqrt {\frac{{{g_i}}}{{{g_j}}}}
+144{r_{ij}}{v_{ij}}\frac{1}{{{g_i}}}) \} }
\\&+ \sum\limits_i {\sum\limits_{\mathop j\limits_{j \ne i} } {{\theta _{ij}}\{  - 24{r_{ii}}\frac{1}{{{g_i}}}\sqrt {\frac{{{g_j}}}{{{g_i}}}}
+ 200{r_{ij}}\frac{1}{{{g_j}}} } }       \\
&  \hspace{80pt} +   \sum\limits_k [{{r_{ik}}{v_{ik}}(120\frac{1}{{\sqrt {{g_i}{g_j}} }}- 144\frac{1}{{{g_i}}}\sqrt {\frac{{{g_j}}}{{{g_i}}}})}
+{{r_{jk}}{v_{ik}}(85\frac{1}{{{g_i}}}+45\frac{1}{{{g_j}}})]\}}
 \\&- 576\sum\limits_i {r_{ii}^3\frac{1}{{{g_i}}}}
  - 576\sum\limits_i {\frac{1}{{{g_i}}}{{\{ \sum\limits_j {{r_{ij}}{v_{ij}}} \} }^3}}
+ \sum\limits_{i,j} {\{ 480r_{ij}^3\frac{1}{{\sqrt {{g_i}{g_j}} }} - 23{r_{ii}}r_{ij}^2\frac{1}{{{g_i}}} - 1728r_{ii}^2{r_{ij}}{v_{ij}}\frac{1}{{{g_i}}}\}}
 \\&+ \sum\limits_{i,j,k} {\{  - 24{r_{ii}}{r_{ik}}{r_{jk}}\frac{1}{{{g_i}}}\sqrt {\frac{{{g_j}}}{{{g_i}}}}
+ 115{r_{ij}}{r_{ik}}{r_{jk}}\frac{1}{{{g_i}}}
+ 1452r_{ik}^2{r_{ij}}{v_{ij}}\frac{1}{{{g_i}}}
  - 1728{r_{ii}}{r_{ij}}{v_{ij}}{r_{ik}}{v_{ik}}\frac{1}{{{g_i}}}\} }
 \\&+ \sum\limits_{i,j,k,l} \{ 120{r_{ik}}{r_{jk}}{r_{il}}{v_{il}}\frac{1}{{\sqrt {{g_i}{g_j}} }}
- 144 \, {r_{ij}}{r_{il}}{r_{jk}}{v_{jk}}\frac{1}{{{g_j}}}\sqrt {\frac{{{g_l}}}{{{g_j}}}}  \\
& \hspace{80pt}  - 40 \, {r_{ik}}{r_{jk}}{r_{il}}{v_{jl}}\frac{1}{{{g_i}}}
+ 720 \, {r_{ij}}{r_{ik}}{v_{ik}}{r_{jl}}{v_{jl}}\frac{1}{{\sqrt {{g_i}{g_j}} }}\} .
      \end{align*}

\section{Formulas of genus-1 3-point functions }
\label{sec:3ptg1}

Explicit formulas for genus-1 3-point functions were used in the proof of Lemma~\ref{lem:thetaijk}
and Lemma~\ref{lem:Omegaij}. Such formulas can be obtained using equations \eqref{eq:g=1 1-pt}, \eqref{eq:ditui} and \eqref{eq:g=1 2-ptij}.

For distinct $i$, $j$, $k$, the genus-1 3-point function $\phi_{ijk}$ has the following form:
\begin{align}
 12{\phi _{ijk}}=&-24{r_{ik}}{r_{jk}}\frac{{{h_j}}}{{{h_i}}}{\phi _i} - 24{r_{ik}}{r_{jk}}\frac{{{h_i}}}{{{h_j}}}{\phi _j} - 24{r_{ij}}{r_{ik}}\frac{{{h_k}}}{{{h_j}}}{\phi _j} - 24{r_{ij}}{r_{jk}}\frac{{{h_k}}}{{{h_i}}}{\phi _i}
- 24{r_{ik}}{r_{ij}}\frac{{{h_j}}}{{{h_k}}}{\phi _k}
              \nonumber\\&- 24{r_{ij}}{r_{jk}}\frac{{{h_i}}}{{{h_k}}}{\phi _k}
- 24{r_{ik}}{r_{jk}}\frac{{{h_j}}}{{{h_k}}}\frac{{{h_i}}}{{{h_k}}}{\phi _k} - 24{r_{ij}}{r_{jk}}\frac{{{h_k}}}{{{h_j}}}\frac{{{h_i}}}{{{h_j}}}{\phi _j} - 24{r_{ik}}{r_{ji}}\frac{{{h_k}}}{{{h_i}}}\frac{{{h_j}}}{{{h_i}}}{\phi _i}
             \nonumber\\& - 24{r_{ij}}\frac{{{h_i}}}{{{h_j}}}{\phi _{jk}} - 24{r_{ik}}\frac{{{h_i}}}{{{h_k}}}{\phi _{jk}} - 24{r_{ij}}\frac{{{h_j}}}{{{h_i}}}{\phi _{ik}} - 24{r_{jk}}\frac{{{h_j}}}{{{h_k}}}{\phi _{ik}}
- 24{r_{ik}}\frac{{{h_k}}}{{{h_i}}}{\phi _{ij}} - 24{r_{jk}}\frac{{{h_k}}}{{{h_j}}}{\phi _{ij}}
              \nonumber
\\&-{\theta _{ik}}{r_{jk}}\frac{{{h_j}}}{{{h_i}}}- \theta _{ki}r_{ij}\frac{{{h_j}}}{{{h_k}}} - {\theta _{jk}}{r_{ik}}\frac{{{h_i}}}{{{h_j}}}
 -\theta _{kj}r_{ij}\frac{{{h_i}}}{{{h_k}}}- {\theta _{ij}}{r_{jk}}\frac{{{h_k}}}{{{h_i}}}  - {\theta _{ji}}{r_{ik}}\frac{{{h_k}}}{h_j}
    \nonumber
\\&+6r_{ij}^2{r_{ik}}\frac{{{h_k}}}{{{h_i}}} + 6r_{ij}^2{r_{jk}}\frac{{{h_k}}}{{{h_j}}} + 6r_{ik}^2{r_{ij}}\frac{{{h_j}}}{{{h_i}}}
 + 6r_{jk}^2{r_{ij}}\frac{{{h_i}}}{{{h_j}}} + 6r_{ik}^2{r_{jk}}\frac{{{h_j}}}{{{h_k}}} + 6r_{jk}^2{r_{ik}}\frac{{{h_i}}}{{{h_k}}}
                             \nonumber\\&+ 12{r_{ij}}{r_{ik}}{r_{jk}}
                            \nonumber
\\&+ \sum\limits_l \{  {r_{il}}{r_{jk}}{r_{kl}}\frac{{{h_i}{h_j}}}{{h_l^2}} + {r_{ik}}{r_{jl}}{r_{kl}}\frac{{{h_i}{h_j}}}{{h_l^2}}
                             + {r_{ik}}{r_{il}}{r_{jl}}\frac{{{h_j}{h_k}}}{{h_l^2}} + {r_{ij}}{r_{il}}{r_{kl}}\frac{{{h_j}{h_k}}}{{h_l^2}} \nonumber
\\&\quad\quad\quad + {r_{ij}}{r_{jl}}{r_{kl}}\frac{{{h_i}{h_k}}}{{h_l^2}} + {r_{il}}{r_{jk}}{r_{jl}}\frac{{{h_i}{h_k}}}{{h_l^2}} - {r_{il}}{r_{jl}}{r_{kl}}\frac{{{h_i}{h_j}{h_k}}}{{h_l^3}}
                            - {r_{ik}}{r_{jk}}{r_{kl}}\frac{{{h_i}{h_j}}}{{{h_k}{h_l}}} \nonumber
\\& \quad\quad\quad- {r_{ij}}{r_{ik}}{r_{il}}\frac{{{h_j}{h_k}}}{{{h_i}{h_l}}} - {r_{ij}}{r_{jk}}{r_{jl}}\frac{{{h_i}{h_k}}}{{{h_j}{h_l}}}
                            - {r_{ij}}{r_{jk}}{r_{kl}}\frac{{{h_i}}}{{{h_l}}} - {r_{ik}}{r_{jk}}{r_{jl}}\frac{{{h_i}}}{{{h_l}}} \nonumber
\\&\quad\quad\quad- {r_{ik}}{r_{jk}}{r_{il}}\frac{{{h_j}}}{{{h_l}}} - {r_{ij}}{r_{ik}}{r_{kl}}\frac{{{h_j}}}{{{h_l}}}
    - {r_{ij}}{r_{il}}{r_{jk}}\frac{{{h_k}}}{{{h_l}}} - {r_{ij}}{r_{ik}}{r_{jl}}\frac{{{h_k}}}{{{h_l}}}\}.
                                \label{g=1 3-ptijk }
                                \end{align}
For $i \neq j$, the genus-1 3-point function $\phi_{iij}$ has the form
\begin{align}
 24{\phi _{iij}}
=&{{\widetilde \Omega }_{ij}} + {\theta _{ij}}\{ 24{r_{ij}} + 4{r_{ii}}\frac{{{h_j}}}{{{h_i}}} - \sum\limits_k {} ({r_{ik}}\frac{{{h_j}}}{{{h_k}}} + {r_{jk}}\frac{{{h_i}}}{{{h_k}}})\}
\nonumber\\& + {\theta _{ji}}(4{r_{ij}}\frac{{h_i^2}}{{h_j^2}}) + \sum\limits_{k \ne i} {{\theta _{ik}}} {r_{jk}}\frac{{{h_j}}}{{{h_i}}} + \sum\limits_{k \ne i} {{\theta _{ki}}} \{ 3{r_{ij}}\frac{{{h_j}{h_k}}}{{h_i^2}} - {r_{ij}}\frac{{{h_j}}}{{{h_k}}} - {r_{jk}}\frac{{{h_i}{h_j}}}{{h_k^2}}\}
\nonumber\\& - 24r_{ii}^2{r_{ij}}\frac{{{h_j}}}{{{h_i}}} - 24r_{ij}^3\frac{{{h_i}}}{{{h_j}}}
\nonumber\\& + \sum\limits_k {} \{ 21{r_{ij}}r_{ik}^2\frac{{{h_j}}}{{{h_i}}} - 2r_{ik}^2{r_{jk}}\frac{{h_i^2{h_j}}}{{h_k^3}} + 4{r_{ii}}{r_{ij}}{r_{ik}}\frac{{{h_j}}}{{{h_k}}} + 2{r_{jk}}r_{ij}^2\frac{{h_i^2}}{{{h_j}{h_k}}}\}
\nonumber\\& - \sum\limits_{k,l} {({r_{ik}}{r_{jl}}{r_{kl}}\frac{{{h_i}{h_j}}}{{h_l^2}} + {r_{ij}}{r_{ik}}{r_{kl}}\frac{{{h_j}}}{{{h_l}}} )}
\nonumber\\& + 96{r_{ii}}{r_{ij}}\frac{{{h_j}}}{{{h_i}}}{\phi _i} -24 {\theta _{ij}}({\phi _i}\frac{{{h_j}}}{{{h_i}}} + {\phi _j}\frac{{{h_i}}}{{{h_j}}}) + 48r_{ij}^2\frac{{h_i^2}}{{h_j^2}}{\phi _j}
 - 24\sum\limits_k {{r_{ij}}{r_{ik}}\frac{{{h_j}}}{{{h_k}}}{\phi _k}}
  \nonumber\\&  + 24\sum\limits_k {{r_{ik}}} \frac{{{h_i}}}{{{h_k}}}{\phi _{jk}}
    + 3{r_{ij}}\frac{{{h_j}}}{h_{i}^{3}} < \tau_{-}^2(S),{\mathcal{E}_{i}} >.
 \label{g=1 3-ptiij }
 \end{align}


\vspace{30pt} \noindent
Xiaobo Liu \\
Beijing International Center for Mathematical Research,\\
Beijing University, Beijing, China.\\
E-mail address: {\it xbliu@math.pku.edu.cn} \\
\& \\
\noindent
Department of Mathematics, \\
University of Notre Dame, \\
Notre Dame,  IN  46556, USA \\
E-mail address: {\it xliu3@nd.edu}\\
\\
\\
Xin Wang \\
School of Mathematical Sciences,  \\
Beijing University, Beijing, China \\
E-mail address:{\it xinwang-1989@163.com}


\begin{thebibliography}{399}

\bibitem[DW]{DW} R. Dijkgraaf and E. Witten,
        {\it Mean field theory, topological field theory, and multimatrix
           models}, Nucl. Phys. B 342 (1990) 486-522.

\bibitem[D]{D}   B. Dubrovin,
        {\it Geometry of 2D topological field theories}, Integrable Systems and Quantum Group,
        Spinger Lecture Notes in Math. 1620 (1996), 120-348.

\bibitem[DLZ]{DLZ}   B. Dubrovin, S. Liu, Y. Zhang,
        {\it On the genus Two Free Energies for Semisimple Frobenius Manifold}, Russian Journal of Mathematical Physics 19 (2012), 273-298.
\bibitem[DZ]{DZ}   B. Dubrovin, Y. Zhang,
        {\it Normal forms of integrable PDEs, Frobenius manifolds and Gromov-Witten invariants}, arxiv:math/0108160.
\bibitem[FJR]{FJR} H.Fan, T.Jarvis, Y.Ruan,
 {\it The Witten equation, mirror symmetry and quantum singularity theory,},
 eprint arxiv:0712.4021, to appear in Ann. of Math.
\bibitem[FLZZ]{FLZZ}   Y. Fu, S. Liu, Y. Zhang, C. Zhou,
        {\it Proof of a Conjecture On the genus Two Free Energies Associated to the $A_{n}$ Singularity},
        arXiv:1205.5990.
\bibitem[Ge]{Ge} E. Getzler,
        {\it Topological recursion relations in genus 2},
       Integrable systems and algebraic geometry (Kobe/kyoto, 1997)
        73-106.
\bibitem[KM]{KM} M. Konsevich, Yu. I. Manin,
        {\it Gromov-Witten classes, quantum cohomology, and enumerative geometry},
       Commun. Math. Phys. 164 (1994), 525-562.

\bibitem[L02]{L02} X. Liu,
    {\it Quantum product on the big phase space and Virasoro conjecture},
    Advances in Mathematics 169 (2002), 313-375.

\bibitem[L06]{L06} X. Liu,
    {\it Idempotents on the big phase space},
    Contemporary Mathematics,vol. 403 (2006), 43-66.

\bibitem[L07]{L07} X. Liu,
    {\it Genus-2 Gromov-Witten invariants for manifolds with semisimple quantum cohomology},
    American Journal of Mathematics, 129 (2007), no.2, 463-498.


\end{thebibliography}
\end{document}